%% file: main.tex
\documentclass[10pt]{amsart}

\usepackage{amssymb}
\usepackage{amsmath}
\input{macros.tex}
\begin{document}

\title[Parapuzzle of the Multibrot set]{Parapuzzle of the Multibrot set\\
  and typical dynamics of unimodal maps}

\author{Artur Avila, Mikhail Lyubich and Weixiao Shen}

\address{
CNRS UMR 7599, Laboratoire de Probabilit\'es et Mod\`eles al\'eatoires\\
Universit\'e Pierre et Marie Curie--Bo\^\i te courrier 188\\
75252--Paris Cedex 05, France
}
\curraddr{IMPA -- Estrada D. Castorina 110, Jardim Bot\^anico, 22460-320
Rio de Janeiro -- Brazil.}
\email{artur@math.sunysb.edu}

\address{Mathematics Department and IMS, Stony Brook University,
Stony Brook, NY 11794, USA}
\email{mlyubich@math.sunysb.edu}
 
\address{Department of Mathematics, University of Toronto, Ontario Canada
M5S 3G3} 
\email{misha@math.toronto.edu}

\address{Mathematics Department,
University of Science and Technology of China, Hefei, 230026,
CHINA} 
\email{wxshen@ustc.edu.cn}

\begin{abstract}

We study the parameter space of unicritical polynomials $f_c:z\mapsto z^d+c$.
For complex parameters, we prove that for Lebesgue almost every $c$, 
the map $f_c$ is either hyperbolic or infinitely renormalizable.  
For real parameters, we prove that for Lebesgue almost every $c$,
the map $f_c$ is either hyperbolic, or Collet-Eckmann, or infinitely renormalizable. 
These results are based on controlling the spacing between
consecutive elements in the ``principal nest'' of parapuzzle pieces.

\end{abstract}

\dedicatory{To memory of Adrien Douady}

\date{\today}

\setcounter{tocdepth}{1}

\maketitle
\thispagestyle{empty}
\input{imsmark}
\SBIMSMark{2008/3}{April 2008}{}

\tableofcontents

\input{intro.tex}
\input{hm.tex}

\input{parapuzzle.tex}

\input{statistics.tex}

\input{outline.tex}

\input{bib.tex}
\end{document}

%% file: macros.tex
\let\newpf\proof \let\proof\relax 
\newenvironment{pf}{\newpf[\proofname]}{\qed\endtrivlist}

\def\SSS{{\mathcal{S}}}

\def\tVV{{\boldsymbol{\VV}}}

\def\ssubset{\Subset}

\def\mo{{\operatorname{m}}}

\def\P{{\mathbb{P}}}

\def\be{\begin{equation}}
\def\ee{\end{equation}}

\def\ba{{\begin{align}}}
\def\ea{{\end{align}}}

\def\n{{\bf n}}
\def\m{{\bf m}}

\def\q{{\bf q}}

\def\dyn{{\mathrm{dyn}}}
\def\p{{\mathrm{par}}}

\def\H{{\mathbb H}}

\def\0{{\mathbf 0}}

\def\cal{\mathcal}

\newtheorem{thm}{Theorem}[section]
\newtheorem{cor}[thm]{Corollary}
\newtheorem{lem}[thm]{Lemma}
\newtheorem{lemma}[thm]{Lemma}

\theoremstyle{remark}
\newtheorem{rem}{Remark}[section]

\numberwithin{equation}{section}

\def \bn {\hfill \\ \smallskip\noindent}

\theoremstyle{definition}

\def\proof{\bn {\bf Proof.} }

\def\note#1
{\marginpar
{\tiny $\leftarrow$
\par
\hfuzz=20pt \hbadness=9000 \hyphenpenalty=-100 \exhyphenpenalty=-100
\pretolerance=-1 \tolerance=9999 \doublehyphendemerits=-100000
\finalhyphendemerits=-100000 \baselineskip=6pt
#1}\hfuzz=1pt}

\newcommand{\di}{\partial}

\newcommand{\ra}{\rightarrow}

\def\ssk{\smallskip}
\def\msk{\medskip}

\def\nin{\noindent}

\def\sm{\setminus}
\newcommand{\diam}{\operatorname{diam}}
\newcommand{\dist}{\operatorname{dist}}

\newcommand{\inter}{\operatorname{int}}
\renewcommand{\mod}{\operatorname{mod}}
\newcommand{\tl}{\tilde}

\newcommand{\supp}{\operatorname{supp}}

\newcommand{\Dil}{\operatorname{Dil}}

\newcommand{\Dom}{\operatorname{Dom}}

\newcommand{\eps}{{\epsilon}}
\newcommand{\De}{{\Delta}}
\newcommand{\de}{{\delta}}
\newcommand{\la}{{\lambda}}
\newcommand{\La}{{\Lambda}}

\newcommand{\Om}{{\Omega}}

\renewcommand{\AA}{{\cal A}}

\newcommand{\CC}{{\cal C}}

\newcommand{\EE}{{\cal E}}

\newcommand{\KK}{{\cal K}}
\newcommand{\LL}{{\cal L}}
\newcommand{\MM}{{\cal M}}

\newcommand{\OO}{{\cal O}}
\newcommand{\PP}{{\cal P}}
\newcommand{\QQ}{{\cal Q}}
\newcommand{\RR}{{\cal R}}

\newcommand{\UU}{{\cal U}}
\newcommand{\VV}{{\cal V}}
\newcommand{\WW}{{\cal W}}
\newcommand{\XX}{{\cal X}}
\newcommand{\YY}{{\cal Y}}

\newcommand{\C}{{\mathbb C}}
\newcommand{\D}{{\mathbb D}}

\newcommand{\R}{{\mathbb R}}
\newcommand{\T}{{\mathbb T}}

\newcommand{\X}{{\mathbb X}}

\def\Bf{{\mathbf{f}}}
\def\Bg{{\mathbf{g}}}

\def\B0{{\mathbf{0}}}
\def\BV{{\mathbf{V}}}
\def\Bv{{\mathbf{v}}}

\catcode`\@=12

\def\Empty{}
\newcommand\oplabel[1]{
  \def\OpArg{#1} \ifx \OpArg\Empty {} \else
  	\label{#1}
  \fi}
		
\newcommand{\comm}[1]{}
\newcommand{\comment}[1]{}

%% file: imsmark.tex
\def\IMSmarkvadjust{0 pt}
\def\IMSmarkhadjust{0 pt}
\def\IMSmarkhpadding{0 pt}
\def\IMSpubltext{Published in modified form:}
\def\SBIMSMark#1#2#3{
 \font\SBF=cmss10 at 10 true pt
 \font\SBI=cmssi10 at 10 true pt
 \setbox0=\hbox{\SBF \hbox to \IMSmarkhpadding{\relax}
                Stony Brook IMS Preprint \##1}
 \setbox2=\hbox to \wd0{\hfil \SBI #2}
 \setbox4=\hbox to \wd0{\hfil \SBI #3}
 \setbox6=\hbox to \wd0{\hss
             \vbox{\hsize=\wd0 \parskip=0pt \baselineskip=10 true pt
                   \copy0 \break%
                   \copy2 \break%
                   \copy4 \break}}
 \dimen0=\ht6   \advance\dimen0 by \vsize \advance\dimen0 by 8 true pt
                \advance\dimen0 by -\pagetotal
	        \advance\dimen0 by \IMSmarkvadjust
 \dimen2=\hsize \advance\dimen2 by .25 true in
	        \advance\dimen2 by \IMSmarkhadjust

%
%
  \openin2=publishd.tex
  \ifeof2\setbox0=\hbox to 0pt{}
  \else 
     \setbox0=\hbox to 3.1 true in{
                \vbox to \ht6{\hsize=3 true in \parskip=0pt  \noindent  
                {\SBI \IMSpubltext}\hfil\break
                \input publishd.tex 
                \vfill}}
  \fi
  \closein2
  \ht0=0pt \dp0=0pt
 \ht6=0pt \dp6=0pt
 \setbox8=\vbox to \dimen0{\vfill \hbox to \dimen2{\copy0 \hss \copy6}}
 \ht8=0pt \dp8=0pt \wd8=0pt
 \copy8
 \message{*** Stony Brook IMS Preprint #1, #2. #3 ***}
}

%% file: intro.tex
\section{Introduction}

In this paper we are concerned with the dynamics of {\it unicritical}
polynomials
\be
f_c : z\mapsto z^d+c,
\ee
where $d \geq 2$, both on the real line (for real values of $c$) and on the
complex plane (in the general case).

Until recently, the dynamical theory of the quadratic family ($d=2$) 
had been developed much deeper than its counterpart for the higher degree 
unicritical polynomials (see \cite{H,M1,L7,S}, \cite{puzzle}--\cite{horseshoe}, \cite{AM1}).  
The reason was that the quadratic maps possess some very special 
geometric features that distinguish them from their higher degree cousins.   
Recently, new tools have been developed  \cite{KL1,KL2,AKLS}
that opened an opportunity  to bring the higher degree case
to the same level of maturity as the quadratic one.%
\footnote{See also \cite{Sm,KSS,BSS} for recent advances in the higher degree case
that use different tools.} 
In this paper that deals with the at most finitely renormalizable case,
combined with forthcoming notes dealing with the
infinitely renormalizable case,
we intend to accomplish this goal.

For $d \geq 2$ fixed,
let $\MM=\MM_d=\{c \in \C,$ the Julia set of $f_c$ is connected$\}$
be the corresponding {\it Multibrot set}.
The dynamics when $c \notin \MM$ is always trivial, so we are
mostly concerned with the description of the dynamics for $c \in \MM$.
When $d$ is odd, the real dynamics is trivial for all $c \in \R$, since
$f_c$ is a homeomorphism, so when discussing real dynamics we will always
assume that $d$ is even.  In this case, for $c \in \MM \cap \R$, $f_c$ is a
unimodal map.

In what follows, various properties of a map $f_c$ will also be attributed to the
corresponding parameter $c$.
For a real $c$, the  map $f_c$ (and the parameter $c$ itself)  are called

\ssk\noindent$\bullet$ {\it regular} if $f_c$ has an attracting periodic orbit;

\ssk\noindent$\bullet$ {\it infinitely renormalizable} if there exist
   periodic intervals of an arbitrarily large (minimal) period;

\ssk\noindent$\bullet$ {\it Collet-Eckmann} if there exist $C>0$ and $\la>1$ such that  
$$  
     |Df^n(c)| \geq C\la^n, \quad n=0,1,2, \dots 
$$ 
Such a map has a unique absolutely continuous invariant measure 
with strong statistical properties. 

\ssk
We can now formulate our main result on the real dynamics: 

\begin{thm} \label {ce}

For almost every $c \in \MM_d \cap \R$, the map $f_c$ is either regular, or
Collet-Eckmann, or infinitely renormalizable.

\end{thm}

\begin{rem}
 In \cite{MN}, Martens and Nowicki described a bigger class of unimodal maps that have an
absolutely continuous invariant measure. 
In \cite {parapuzzle}, it was proved that for almost every real $c\in \MM_2$,
the quadratic polynomial $f_c$  is either regular, or Martens-Nowicki, or infinitely renormalizable.
The Martens-Nowicki property was then replaced in \cite {AM1} with the much stronger Collet-Eckmann property, 
thus providing us with Theorem \ref{ce} in the quadratic case.
\end{rem}

\begin{rem}

With Theorem \ref{ce} in hands, 
we can go further, in the same way as in the quadratic case, to show that  
the whole fine statistical description of the dynamics of real quadratic maps 
\cite {AM1}, \cite {AM4} is valid in  the  higher degree case  as well.

\end{rem}

\begin{rem}
In the forthcoming notes, the above result will be complemented by showing that the set of
infinitely renormalizable parameters in $\MM_d\cap \R$  has zero Lebesgue measure.
(In the quadratic case, this was proved in \cite{horseshoe}.)
\end{rem}

To stress the difference between the quadratic and the higher degree cases, let
us mention one consequence of Theorem \ref{ce}.
Recall that a {\it wild attractor} for a unimodal map is a measure-theoretic
attractor (in the sense of Milnor \cite{M-attractors}) which is not a topological attractor. 
There are no wild attractors in the quadratic family \cite {attractors},
but they do exist for a sufficiently high even criticality $d$ \cite {BKNS}. Moreover,
if $d$ is big enough, the set of parameters $c\in \MM_d\cap \R$ for which 
the wild attractor exists contains a Cantor set.

\begin{cor}
For any even criticality $d$, 
the set of parameters $c \in \MM_d \cap \R$ for which the wild attractor exists
has zero Lebesgue measure.

\end{cor}

\begin{rem}

One can show (using the estimates of this paper)
that for large even $d$, the
set of parameters $c\in \MM_d\cap \R$ for which the wild attractor exists has positive Hausdorff
dimension.

\end{rem}

\begin{rem}

In \cite {BSS}, it was proved that for almost every
$c \in \MM_d \cap \R$, the map $f_c$ admits a physical measure%
\footnote{A  measure $\mu$ is called {\it physical} if 
 the Birkhoff averages of Lebesgue  almost all orbits converge to $\mu$.}, 
which is either
supported on an attracting periodic orbit, 
or is absolutely continuous (but $f_c$ is not necessarily Collet-Eckmann), 
or is supported on a uniquely ergodic Cantor set coinciding with the postcritical
set (this possibility contains strictly the infinitely renormalizable case
and the case of wild attractors).

\end{rem}

The set of non-regular, non-infinitely renormalizable real parameters does
have positive Lebesgue measure \cite {J}, \cite {BC}.  The situation is
quite different for complex parameters:

\begin{thm} \label {measure zero}

For almost any $c \in \C$, the map $f_c: z\mapsto z^d+c$
is either hyperbolic or infinitely renormalizable.

\end{thm}

This result was proved for quadratic maps by Shishikura, see a sketch in
\cite {S} (see also \cite {AM3} for a proof closer to this paper).  We
actually prove the following estimate:

\begin{thm} \label {refined}

Let $f_c$ be a non-renormalizable map with all fixed points
repelling.\footnote {The result still holds under the assumption that $f_c$
is not infinitely renormalizable, and has all periodic orbits repelling.
The argument for this generalization (which is more subtle than the usual
application of the renormalization operator) is indicated in
Remark \ref {gen}.}
Then $c$ is not a density point of $\MM_d$.

\end{thm}

Again, the last two results are more surprising in the case $d>2$: 
while all the quadratic maps in question 
have the Julia set of zero  area (\cite {L7}, \cite {S}), 
it is conceivable that there exist higher degree non-renormalizable unicritical polynomials
with the Julia set of positive measure.
So, in the quadratic case the phase-parameter dictionary works in the natural way:
zero area of Julia sets of the class of maps under consideration
translates into zero area of the corresponding set of parameters.
On the other hand, in the higher degree case, the phase-parameter correspondence is
more subtle.  

\ssk
An important special feature of the quadratic maps essentially exploited in
the previous work  
is the {\it decay of geometry of the principal nest}, see \cite{puzzle,parapuzzle,AM1}.
In this paper we demonstrate that, though in the higher degree case this property fails in general,
it is satisfied for {\it almost all} non-regular non-infinitely renormalizable parameters.  
This is the key to all of the above results.

\comm{
The Multibrot set contains infinitely many copies of itself,
each arising from some renormalization combinatorics, all of them
canonically homeomorphic (via the {\it straightening}) to $\MM$.
Let us call {\it maximal} a copy of the Multibrot set which is not
contained properly in any other copy.  A maximal copy is called
{\it satellite} if it is attached to the {\it main cardioid}, 
the set of parameters with a neutral fixed point, 
and is called {\it primitive} otherwise.  The straightening is
quasiconformal for primitive copies, and for satellite copies it is
quasiconformal outside any neighborhood of the {\it root} which attaches it
to the main cardioid.  Let us say that a Multibrot copy has a $K$-bounded
shape if the straightening admits a $K$-qc extension to the whole complex
plane.

For several purposes, it is convenient to organize
primitive copies in families with certain common combinatorial features. 
Families of copies belonging to the same {\it Misiurewicz wake} were
considered in \cite {parapuzzle} (in the case $d=2$).  Consideration of
Misiurewicz wakes is certainly important for real dynamics: all primitive
copies intersecting the real line belong to the same Misiurewicz wake.

\begin{thm}

All primitive Multibrot copies in a Misiurewicz wake $\OO$
have $K$-bounded shape, where $K=K(\OO)$.

\end{thm}

This result was proved for the quadratic family in \cite {parapuzzle}.
The particular case of the real Misiurewicz wake plays
an important role in the analysis of the real renormalization horseshoe in
\cite {ALS}.
} 

{\bf Acknowledgment:}  We thank all the institutions and foundations
that have supported us in the course of this work: 
Simons Mathematics and Physics  Endowment,  Fields Institute, NSF, NSERC,
University of Toronto, Warwick Mathematics Institute.  W.S. acknowledges
support by the SRFDP grant No. 20070358058 and the 973 program grant No.
2006CB805900.  This research was
partially conducted during the period A.A. served as a Clay Research Fellow.

\subsection{Notations}

\nin
$\D_r= \{ z:\,  |z|<r\}$, $\D=\D_1$, $\T=\partial \D$.

If $S$ is a hyperbolic Riemann surface, let $\dist_S$ be the hyperbolic
metric in $S$ (with the usual normalization, so that in the upper-half
plane $\H$ we have $\dist_\H(i,ai)=|\log a|$, $a>0$).  The diameter of a
subset $X \subset S$ with respect to $\dist_S$ will be denoted $\diam_S X$.

\nin
$K(f)$ is the {\it filled Julia set} of $f$.\\ 
$J(f)=\di K(f)$ is its {\it Julia set}.

\nin
$\Dil(h)$ stands for the {\it dilatation} of a quasiconformal map $h$.

\nin
$\mod(A)$ stands for the modulus of the annulus $A$.

\nin
{\it Pullbacks} of an {\it open} topological disk $V$ under $f$ are
connected components of $f^{-1}(V)$.

\nin
{\it Pullbacks} of a {\it closed} disk $V$ are closures of
pullbacks of $\inter V$.

%% file: hm.tex
\section{Holomorphic motions and a phase-parameter lemma}

Let $\Lambda \subset \C$ be a Jordan disk.
A {\it holomorphic motion} over $\Lambda$ (with base point
$\lambda_0 \in \Lambda$) of some set $Z \subset \C$
is a family of injective maps $h_\lambda:Z \to \C$, $\lambda \in \Lambda$,
such that for every $z \in Z$, the ``trajectory'' (or the ``orbit'')
$\lambda \mapsto h_\lambda(z)$ 
is holomorphic in $\la$ and $h_{\lambda_0}(z)=z$.
Given such a holomorphic motion, we let $Z_\la=h_\la(Z)$. 

The central result in the theory of holomorphic motions  is the {\it $\la$-lemma}. 
It consists of two parts: extension and quasiconformality. 
The Extension Theorem (in its strongest version which is due to Slodkowski \cite {Sl})
says that a holomorphic motion $h_\lambda:Z \to \C$ over a Jordan disk $\Lambda$ 
can be always extended to a holomorphic motion $\hat h_\lambda:\C \to \C$ of the whole plane 
over the same $\Lambda$. 
The Quasiconformality Theorem (Ma\~n\'e-Sad-Sullivan \cite{MSS}) 
states that each $\hat h_\lambda$ is quasiconformal and
$$
   \log \Dil(h_\lambda) \leq \dist_\Lambda (\lambda_0,\lambda).
$$

We say that a holomorphic motion $h_\la: Z\ra \C$ is {\it continuous up to the boundary}
if the map $(\la, z)\mapsto h_\la(z)$ extends continuously to $\bar \La\times Z$. 
A holomorphic motion $h_\la$ of a Jordan curve $T$ over $\La$ which is continuous up to the boundary
will be called a {\it tubing} of $T$ over $\La$.
Under these circumstances, a {\it diagonal} to the tubing is a holomorphic function
$\psi$ in a neighborhood of  $\bar \Lambda$  
satisfying the following properties:

\ssk\nin (D1)
For $\lambda \in \Lambda$, $\psi(\lambda)$ belongs to the bounded component
of $\C \setminus T_\lambda$, and for $\lambda \in \partial \Lambda$,
$\psi(\la)\in T_\lambda$.

\ssk\nin  (D2)
   For any  $\la\in \di \La$, the point $\psi(\la)$ has only one preimage $\gamma(\la) \in T$ under $h_\la|T$;

\ssk\nin (D3)
 The holomorphic motion of a neighborhood of $\gamma(\la)$ in $T$ admits an extension over
a neigborhood of $\la$;

\ssk\nin (D4) The graph of $\psi$  crosses the orbit of $\gamma(\la)$ transversally at $\psi(\la)$; 

\ssk\nin  (D5)
The map $\gamma: \partial \Lambda \to T$ has degree $1$.  

\begin{rem}
   Note that properties (D3) and (D4) imply that $\gamma: \di\La \ra T$ is continuous, 
    so that, (D5) makes sense.
\end{rem}

Given a set $Z$ contained in the closed Jordan disk bounded by $T$, 
 we say that a  holomorphic (and continuous up to the boundary) motion $H_\la$  of $Z$ over $\La$
{\it fits to the tubing} of $T$ if for every $\lambda \in \overline \Lambda$, 
we have $H_\lambda(z)=h_\lambda(z)$ for $z \in Z \cap T$, while 
$H_\lambda(z) \notin h_\lambda(T)$ for $z \in Z \setminus T$.


\begin{lem}\label{Arg Prin}
Let $h_\la: Z\ra \C$ be a holomorphic motion over a Jordan disk $\La$  continuous up to the boundary
that fits to a tubing of a Jordan curve $T$.  Let $\psi$ be a diagonal to
this tubing.
Then for each point $z\in Z$ there exists a unique parameter $\la=\chi(z) \in \Lambda$
such that $h_\la (z)=\psi(\la)$.  The map $\chi: Z\mapsto \bar \La$
is continuous and injective.
Moreover, if $z \in \inter Z$ and $h_{\chi(z)}$ is
locally $K$-quasiconformal at $z$ then $\chi$ is locally
$K$-quasiconformal at $z$.
\end{lem}

\begin{pf}
  Without loss of generality, we can assume that $\La$ is the unit disk with the base point at the origin.
By assumptions (D3)-(D4), $\gamma: \di\La\ra T$ is a local homeomorphism. By (D5), it has degree 1, so that
it is a homeomorphism. Letting $\chi|(Z \cap T)=\gamma^{-1}$, we see that the
first assertion is valid for $z\in Z \cap T$.   

Let $z\in Z\sm T$. By applying an appropriate family of affine changes of variable,
we can be reduced to the case when  $h_\la(z)=0$, $\la\in \bar\D$.   
 
 Let us consider a torus $\T_0^2 = \T\times T\subset \T\times \C$.
Let us deform  it in $\T\times \C$ as follows:
$$
        H_r : \T_0^2\ra \T_r^2,\  (\la, z) \ra (\la, h_{r\la}(z)),\quad \la\in \T, \ z\in T, \ 0\leq r\leq 1.  
$$   
Since the origin fits to the tubing of $T$,
the deformations never cross the core circle $\T\times \{0\}$.

Let us consider a family of curves $\psi_r: \T\ra \C $,
$\psi_r(\la) = h_{r\la}\circ \gamma(\la)$,  $0\leq r\leq 1$.
Note that the graph of $\psi_r$ is a curve in $\T^2_r$ obtained by applying 
the homotopy $H_r$ to the graph of $\psi_0$. 
Since the $\T^2_r$ are disjoint from the core circle, 
the curves $\psi_r$ never pass through the origin and hence have  the same winding number around it.  
Since $\psi_0=\gamma$, by (D5) this winding number is equal to 1. 
But $\psi_1=\psi|\, \T$ by Definition of $\gamma$ (D2). By the Argument Principle,  
$\psi$ has a single root in $\D$, which proves the first assertion.

Any point $\la\in \La$ has at most one preimage under $\chi$ since the maps $h_\la$ are injections.
A point $\la\in \di \La$ has only one preimage $\chi^{-1}(\la)=\gamma(\la)$  by (D2) and the assumption
that the motion of $Z$ fits to the tubing of $T$.  The graph of $\chi$ is
the set of solutions $(z,\lambda)$ of $h_\lambda(z)=\psi(\lambda)$, which is
clearly closed in $Z \times \overline \Lambda$, so $\chi$ is
continuous.

Local quasiconformality of $\chi$ follows from the $\la$-lemma (see Corollary 2.1 of \cite{parapuzzle}).  
\end{pf}

We will often encounter the situation when  $Z$ contains an annulus $A$, 
and we want to obtain a lower bound on $\mod(\chi(A))$. A trivial  bound
$$
  \mod(\chi(A)) \geq K^{-1} \mod(A),\quad \mathrm {where}\quad  K=\sup_{\lambda \in \Lambda} \Dil(h_\lambda|A),
$$
will sometimes be sufficient.  However, since the dilatation of the holomorphic motion can blow up
as $\la\to \di  \La$, it will not cover all of our needs.
Then we will make use of the following generalization of Corollary 4.5 of \cite {parapuzzle}.

\begin{lemma} \label {1}

Under the above circumstances, let $X \subset Z \setminus T$, and let
$U_\lambda$ be the bounded component of $\C \setminus T_\lambda$. Then:
\begin{enumerate}

\item 
 There exists a $\de_0>0$ such that 
if $\diam_{U_\la} X_\la <\delta\leq \de_0$ for every $\lambda \in \Lambda$ 
then $\diam_\Lambda \chi(X)<\epsilon(\de)$,\footnote {On the other hand, one
can show that the statement is false for large $\delta$.}
where $\eps(\de)\to 0$ as $\de\to 0$.

\item Assume that $X$ is {\it connected} and 
   $\diam_U X \leq M$.
Assume also that for some $K>1$ and for every
$\lambda \in \Lambda$, the map $h_\lambda:X \to \C$
extends to a $K$-qc homeomorphism
$U \to U_\lambda$.
Then  $\diam_\La \chi(X)\leq C=C(M,K)$.
\end{enumerate}

\end{lemma}

\begin{pf}

{\it First statement}.  It is enough to consider the case where $X$
consists of two points.  By a holomorphic change of
coordinates $(\lambda,z) \mapsto (\lambda,\phi_\lambda(z))$ where
$\phi_\lambda:\C \to \C$ is affine, we may assume that $X_\la =\{0,1\}$ 
for all $\lambda \in \Lambda$. 
If $\diam_{U_\la} X_\la $ is small, then $\D_{2R} \subset U_{\la}$ 
for some large $R>1$ and all $\la\in \La$.
Let $\tilde h$ be the holomorphic motion of $T \cup \D_R $ obtained
by setting $\tilde h_\lambda(z)=h_\lambda(z)$ for $z \in T$ and
$\tilde h_\lambda(z)=z$ for $z \in \D_R $.  Notice that $\tilde h$ also fits
to the tubing and $\tl \chi(X)= \chi (X)$.  
  
Since $\tilde h$ is holomorphic at $\D_R$,
$\tilde \chi= \psi^{-1}:\tilde W \to \Lambda$ is also holomorphic on $\D_R$.
Hence $\mod(\tilde \chi(\D_R\sm \D ))=\displaystyle{ \frac{1}{2\pi}\log R}$ and   
$$
    \diam_\La (\chi(X)) \leq \diam_{\tl\chi (\D_R) } \tilde\chi(X) = O(1/R).
$$

\msk
{\it Second statement}.  We will use the uniform equicontinuity of $K$-qc maps
with respect to the hyperbolic metric:
For any $K$-qc map $\phi:S \to \tilde S$ between hyperbolic Riemann surfaces,
$$
  \dist(x,y)<\eta \implies \dist(\phi(x), \phi(y)) < \de (K,\eta),
$$
where $\de(K,\eta)\to 0$ as $\eta \to 0$. 
   
Let us select $\eta=\eta(K)$ so that $\de(K,\eta) < \de_0$,
where $\de_0$ comes from the first statement. 
We can cover $X$ by $N=N(\eta,M)$ sets $X_1$,...,$X_N$ of
hyperbolic diameter in $U$  bounded by $\eta$. 
Then the first statement is applicable to each $X_i$,
so that $\diam_\La \chi (X_i) <\eps_0= \eps(\de_0)$.
Since $X$ is connected,  
$$
  \diam_\La \chi(X)\leq  \sum \diam_\La \chi(X_i) < N \eps_0,
$$
and we are done.
\end{pf}

We will need one lemma on lifting of a holomorphic motion by a family of branched coverings. 

\begin{lem}\label{lifting of motions}
  Let $h_\la: Z\ra Z_\la$ be a holomorphic motion over a pointed disk $(\La, \la_0)$,
and let $f_\la: U_\la' \ra U_\la$ be a holomorphic family of branched coverings of degree $d$
such that $U_\la\supset Z_\la$.
\footnote{We assume (as part of the definition of a ``holomorphic family'') 
   that $\cup U_\la$ and $\cup U_\la'$ are open subsets of $\C^2$.}
Let $\Om\Subset \La$ be an open Jordan disk containing $\la_0$  such that for $\la\in \Om$, 
the sets $Z_\la$  do not contain the critical values of $f_\la$.   
Then $h_\la$ over $\Om$ lifts by $f_\la$ to a  holomorphic motion $h_\la'$ 
continuous up to the boundary.
\end{lem}

\begin{pf}
   Each orbit $Z(z)=\{(\lambda,h_\la(z)): \ \la\in \La\}\subset \C^2 $
of the motion $h_\la$ lifts to a variety 
$X(z)  = \{ (\la,z):\ (\lambda,f_\la(z))\in Z(z) \}$ which properly projects to $\La$ with degree $d$. 
Since for $\la\in \Omega$,
$Z_\la$ do not contain critical values of $f_\la$, these varieties are unbranched over $\Om$ and
hence form a holomorphic motion $h_\la': Y\ra Y_\la$ over it. 
All we need to show is that it is continuous up to the boundary of $\Om$. 

It is enough to show that for any compact $K \subset Y$,
the family $\{\lambda \mapsto
h_\lambda'(y)\}_{y \in K}$ is uniformly equicontinuous over $\Omega$.
Let $y_n \in K$, $\sigma_n \subset
\Omega$, $\sigma_n$ an arc of diameter at most $1/n$, and let
$B_n=\{h_\lambda'(y_n):\ \lambda \in \sigma_n\}$.  We must show that the diameter
of $B_n$ shrinks to $0$.  We may assume that $y_n \to y \in K$ and
$\sigma_n \to \lambda \in \overline \Omega$ in the Hausdorff topology.
Then, for any $\epsilon>0$, for large $n$,
$B_n$ lies within an $\epsilon$-neighborhood of
$f_\lambda^{-1}(h_\lambda(f_{\lambda_0}(y)))$.  Since
$f_\lambda^{-1}(h_\lambda(f_{\lambda_0}(y)))$ has at most $d$ elements and
$B_n$ is connected, this implies that $B_n$ has diameter at most
$2 d \epsilon$, as desired.
\comm{
Let us consider a trajectory $Y(\zeta)$, $\zeta\in X (z)$, of the motion $h_\la'$, 
where  $z=f_{\la_0}(\zeta)$. 
Since the projections $X(z)\ra \Om$ are branched coverings of degree $d$,
the extremal distance from $Y(\zeta)$ to the boundary of $X(z)$ is bounded away from $0$. 
Hence each disk $Y(\zeta)$ is contained in some disk $\hat Y(\zeta)\subset X(z)$ with 
\be\label{eps}
\mod(\hat Y(\zeta)\sm Y(\zeta))\geq \eps >0.
\ee
Let us uniformize the disks $\hat Y(\zeta)$ by the unit disk,
$\phi_\zeta: (\D,0)\ra (\hat Y(\zeta), b(\zeta))$, 
where the base point $b(\zeta)\in Y(\zeta)$ projects to $\la_0$.
We obtain a normal family of maps $\phi_\zeta: \D\ra \C^2$. 
Moreover, by (\ref{eps}),
the disks $D_\zeta = \phi_\zeta^{-1}(Y(\zeta))$ are contained in some disk $\D_r$, $r<1$. 
Thus, the family of maps  $\phi_\zeta: \D\ra \C^2$ is equicontinuous. 

It easily follows from the local structure of analytic varieties that
the closures of $Y(\zeta)$ are graphs of some continuous functions
$\psi_\zeta: \bar \Om\ra \C$.
Moreover, equicontinuity of the family of maps $\phi_\zeta$ implies that
the family of functions $\psi_\zeta$ is also equicontinuous. 
Hence the map $(\la, \zeta) \mapsto h_\la'(\zeta ) = \psi_\zeta(\la)$ 
is continuous on  $\bar \Om\times Y$.  
}
\end{pf}

%% file: parapuzzle.tex
\section{Puzzle and parapuzzle}

\subsection{Parameter and dynamical B\"ottcher coordinates}
 The basic dynamical theory of the unicritical family $z\mapsto z^d+c$ (see \cite{Sch1}) 
is similar to the basic theory of the quadratic family  (see \cite{DH,M2}).
For further reference, we recall here the main objects of the theory and set up notations.

The {\it (dynamical) B\"ottcher function} $B_c$ conformally conjugates $f_c$ near $\infty$ to
$z\mapsto z^d$. The {\it Green function} $G_c=\log |B_c|$ extends harmonically to 
$\C\sm K(f_c)$. Its level sets $\{ G_c=\xi \}$
are called {\it (dynamical) equipotentials} $\EE^\dyn_\xi= \EE^\dyn_\xi(c)$. They form an invariant foliation
with singularities at the precritical points
(at each singularity, the equipotential looks locally like the intersection of $d$ lines). 
 Let 
$$
   \De_c=\{z: G_c(z)>G_c(0)\}.
$$ 
It is  the maximal neigborhood of $\infty$ saturated by the equipotentials
on which the foliation is non-singular.  

The gradient lines of $G_c$ coming from infinity are called {\it (dynamical) external rays}.
They form a foliation of $\C \setminus K(f)$ slit along the gradient lines emerging from the critical
points of the Green function. 
The argument (``angle'') of $B_c$ is constant on each ray. The ray of angle $\theta$ is denoted as
$\RR^\dyn_\theta=\RR^\dyn_\theta(c)$.

If the Julia set of $f$ is connected, the B\"ottcher function extends analytically to the whole
basin of infinity, $\C\sm K(f_c)$, and maps it conformally onto $\C\sm
\overline \D$.  

Otherwise, $B_c$ extends analytically to the domain $\De_c$,
and maps it conformally onto $\C\sm \D_{\rho(c)}$, where $\rho(c)= e^{G_c(0)}>1$.
In this case, the function 
\be\label{basic}
        B_\MM(c)= B_c(c)
\ee
is well defined and is called the {\it (parameter) B\"ottcher function}.
It provides us with the Riemann mapping $\C\sm \MM\ra \C\sm \overline\D$.
This {\it basic relation} between the dynamical and parameter B\"ottcher coordinates/Riemann mappings 
is the foundation of the phase-parameter correspondences for the unicritical families of polynomials. 

The {\it (parameter) equipotentials and external rays}, $\EE^\p_\xi$ and $\RR^\p_\theta$, 
 are the level sets and the gradient lines 
of the parameter Green function $G_M(c)=\log |B_\MM(c)|$. They form two (non-singular) orthogonal foliation
on $\C\sm \MM$. By basic relation (\ref{basic}),
\begin{itemize}
  \item  $c\in \EE^\p_\xi$ iff   $c\in \EE^\dyn_\xi(c)$;
  \item  $c\in \RR^\p_\theta$ iff $c\in \RR^\dyn_\theta(c)$.
\end{itemize}
(In each line, the first and the last ``$c$'' stand for the {\it parameter},
while the intermediate one stands for the {\it critical value}.)  

\ssk
Let 
$${\mathbf F}= \{(c, z)\in \C^2: \ z\in \C\sm K(f_c)\} ;\quad 
      {\boldsymbol \De}= \{(c, z)\in \C^2: \ z\in \De_c\};
$$
these are open sets in $\C^2$.
Let us also consider the critical set 
$$
     \CC_- = \{ (c,z)\in {\mathbf F} :\ \exists n\geq 0, \ f_c^n(z)=0\};
$$
it is an analytic subvariety in ${\mathbf F}$.  
The B\"ottcher function
$$
{\mathbf B}: {\boldsymbol \De}\to \C\sm \overline\D, \quad (c,z)\mapsto B_c(z),
$$  
is a local holomorphic submersion, so that, 
its level sets form a holomorphic foliation of
${\boldsymbol \De}$. Moreover, this foliation is transverse to the vertical foliation of $\C^2$,
and thus determines a local holomorphic motion near any point $(c,z)\in {\boldsymbol \De}$. 

Pulling this foliation back by the fiberwise dynamics $\Bf: (c,z)\mapsto (c,f_c z)$,
we obtain a holomorphic foliation on ${\mathbf  F}$ with singularities on $\CC_-$.
It determines a local holomorphic motion near any point $(c,z)\in {\mathbf F}\sm \CC_-$
that we call the {\it B\"ottcher motion}.

We say that some holomorphic motion over parameter domain 
{\it matches the B\"ottcher motion} or {\it respects the B\"ottcher
coordinate}
if on the basin of infinity it coincides with the B\"ottcher motion. 
Such a motion preserves the external angles and heights of the points in the
basin of 
infinity. 

\subsection{Transversality to the diagonal} \label {transverse}

Let $\CC_1=\{(c,c): \ c\in \C\sm M\}$.

\begin{lem}\label{B motion}
   Near any point $(c,c)\in \CC_1$, the B\"ottcher motion is well defined
and is transverse to $\CC_1$.  
\end{lem} 

\begin{pf}
  The B\"ottcher motion is well defined since $\CC_1\cap \CC_-=\emptyset$.
It is transverse to $\CC_1$ since the B\"ottcher funcion ${\mathbf B}|\CC_1$
is non-singular
(as it conformally maps $\CC_1$ onto $\C\sm M$). 
\end{pf}

Let $c_0$ be a Misiurewicz parameter, i.e., there is a repelling periodic
point $a_0$, of period $q$,
such that
$f_{c_0}^n(0)=a_0$ for some $n \geq 1$, assumed to be minimal with this
property.  There are finitely
many (and at least $2$) dynamical rays $\RR^\dyn_{\theta_i}(c_0)$ landing at
$c_0$.
Through a neighborhood of $c_0$, the B\"ottcher motion of these dynamical
rays is well defined, see
Lemma B.1 of \cite {GM} and Lemma 2.2 of \cite {Sch1}.  Their common landing
point $p(c)$
is just the analytic continuation of $c_0$ as a preperiodic point (that is,
$p(c)$ is
the solution of $f_c^{n-1}(z)=f_c^{n-1+q}(z)$ near $c_0$).

\comm{
On the other hand, the periodic point together
with its preimages $f_{c_0}^{-m}(a_0)$, $1 \leq m \leq n-1$, move
holomorphically near $c_0$.  This provides us
with the extension of the B\"ottcher motion of the dynamical rays
$\RR^\dyn_{\theta_i}(c)$
to the closure, which is obtained by adding the common landing point $p(c)$
of the
$\RR^\dyn_{\theta_i}(c)$.
}

\begin{lem}[compare \cite{vS}]\label{Misiurewicz}
   Let $c_0$ be a Misiurewicz parameter as above.
Then the curve $c\mapsto p(c)$ is transverse to the diagonal
$c\mapsto c$ at $c_0$.   
\end{lem}

\begin{pf}
Let us consider one of the dynamical rays $\RR^\dyn_{\theta_i}(c)$ landing
at $p(c)$
which moves holomorphically under the B\"ottcher motion $h_c=B^{-1}_c \circ
B_{c_0}$.
By the basic dynamical-parameter relation, 
$\{c: \ c \in \RR^\dyn_{\theta_i} (c)\}$ is a parameter ray
$\RR^\p_{\theta_i}$
landing at $c$. Moreover, the map $\gamma: c\mapsto h_c^{-1}(c)$ is a
homeomorphism from 
 $\RR^\p_{\theta_i}$ to $\RR^\dyn_{\theta_i} (c_0)$.  But if the curves
 $c\mapsto p(c)$ and $c\mapsto c$ had tangency of order $d\geq 1$ at $c_0$ 
 then each point on  $\RR^\dyn_\theta (c_0)$ would have $d+1$ preimages
 under $\gamma$
(compare Lemma 9.1 of \cite {ALM}) -- contradiction.  
\end{pf}

\subsection{Parabolic wakes}

Let $\AA$ stand for the set of parameters $c$ for which the map $f_c:
z\mapsto z^d+c$ has an attracting fixed point $\alpha_c$. In the quadratic case,
it is a domain bounded by the main cardioid of the Mandelbrot set.
In the higher degree case, $\AA$ is a domain bounded by a Jordan
curve with $d-1$ cusps. 

The set $\MM \setminus \overline \AA$ is disconnected.  
The closures of the connected components of $\MM \setminus \overline \AA$ 
are called {\it (parabolic) limbs} of $M$. 
Each limb $\LL$ intersects $\overline \AA$ at a single point called the
{\it root} $r=r_\LL$ of the limb.  
The map $f_r$ has a parabolic fixed point with
some multiplier $e^{2 \pi i {\bf p}/\q}$.
There are two parameter external rays landing at the root. 
Their union with $r$ divides $\C$ into two (open) connected components: the one
containing $\LL\sm \{r\}$ is called a {\it parabolic wake} $\WW=\WW_\LL$
(see \cite{DH,M2,Sch1}).

For $c\in \LL$,  the map $f_c$ has a unique 
{\it dividing} fixed point $\alpha_c$.
There are $\q$ external rays  $\RR_i^\dyn(c)$ landing at this point
which are cyclically permuted by $f_c$ with {\it combinatorial rotation number} ${\bf p}/\q$. 
This configuration of $\q$ rays, together with the $\alpha$-fixed point, 
moves holomorphically over the whole parabolic wake $\WW$.
We let 
$$
   \Gamma^0=\Gamma^0(c)= \cup \RR^\dyn_i(c).
$$ 

Given some height $\xi>0$, let $\WW(\xi)$ stand for the domain obtained by truncating the
parabolic wake $\WW$
by the parameter equipotential $\EE^\p_{\xi}$ of height $\xi$.
For $c\in \WW(d\xi)$, the {\it Yoccoz puzzle pieces $Y^0_j=Y^0_j(c)$ of depth} $0$ 
are obtained by taking the closure of the
connected components of $\displaystyle{\C \setminus (\Gamma^0 (c) \cup \{\alpha_c\})}$
truncated by the dynamical equipotential  $\EE^\dyn_\xi(c)$ of height $\xi$ 
(where the piece containing $0$ is also denoted $Y^0$). 
This configuration of $\q$ puzzle pieces moves holomorphically over $\WW(d\xi)$.

Since $f_c(0)\not\in \Gamma^0(c) \cup \EE^\dyn_\xi(c)$ for $c\in \WW(\xi)$,
the $f_c$-preimages of the rays $\Gamma^0(c)$ move holomorphically over $\WW(\xi)$,
and so do $f_c^{-1}(\EE^\dyn_\xi(c))=\EE^\dyn_{\xi/d}$ and the $f_c$-preimages of $\alpha_c$.  
The closures of the components of $\C\sm f_c^{-1}(\Gamma^0(c)\cup \alpha_c)$
truncated by the equipotential $\EE^\dyn_{\xi/d}(c)$ are called Yoccoz puzzle pieces of depth 1,
and are denoted $Y^1_j$ (where the one containing $0$ is also denoted $Y^1$). 

We now fix some height $\xi$ (say, $\xi=1$): the moduli bounds in what
follows will depend on this choice, but it will not be explicitly indicated.

\subsection{Satellite copies of $\MM$}

Let  
$$
   \MM_\LL = \{ c\in \LL:\; f^{\q n}(0)\in Y^1,\quad n=0,1,2,\dots \}.
$$ 
This set is canonically homeomorphic to $\MM$, and is called a {\it satellite copy} 
of the Multibrot set (see \cite{DH-pl,Sch2}).
The maps $f_c$ with $c\in \MM_\LL\sm \{r_\LL\}$ (and the corresponding parameters $c$)
are called {\it satellite renormalizable}.

\subsection{Decorations and Misiurewicz wakes}\label{Mis}

Removing the satellite copy from the limb $\LL$ disconnects it 
into countably many components,
each attached to $\MM_\LL$ at a Misiurewicz parameter $c_*$
such that $f_{c_*}^{\n \q}(0) \in f_{c_*}^{-1}(\alpha) \setminus
\{\alpha\}$ for some $\n >0$.  
The closures of these sets are called {\it decorations}.

There are $\q$ rays landing at $c_*$, dividing
$\C$ into $\q-1$ {\it Misiurewicz wakes} and the component containing $\MM_\LL\sm \{c_*\}$. 
The above number $\n$ is called the {\it level} 
of the Misiurewicz wake and the corresponding decoration.   

For $c$ in the Misiurewicz wake,
the level $\n$ is determined as the minimal natural number $n$ such that
$f_c^{n\q}(0)$ belongs to some $Y^1_j(c) \neq Y^1(c)$.
Let $\OO^\n_k$        
stand for the Misiurewicz wakes 
truncated by the parameter equipotential of height  $\xi/d^{\q\n}$.  
Obviously, truncated Misiurewicz wakes are
compactly contained in the correponding truncated parabolic wake $\WW(d\xi)$,
and the Misiurewicz wakes attached to different roots have disjoint closures.

Define
$$
  \Omega^\n =   \WW (\xi / d^{\n\q-1} ) \sm \bigcup_{\m< \n }\bigcup_k
\overline\OO^\m_k . 
$$
It is an open Jordan disk containing $\overline \OO^\n_k$.

For the further understanding of the wakes, we need to go deeper into the puzzle.

The {\it Yoccoz puzzle pieces of depth} $n$  are the
pullbacks of Yoccoz puzzle pieces of depth $0$ under $f^n$.
The puzzle pieces of depth $n$ will be denoted by $Y^n_j$, where the
labels $j$ stand for the angles of the external rays that bound $Y^n_j$.
They form a tiling of the neighborhood of $K(f)$ bounded by the equipotential 
of height $\xi /d^n$. 

We also let $Y^n$ stand for the {\it critical} puzzle piece of depth $n$, i.e., $Y^n\ni 0$,
while $Y^n_v$ stand for the puzzle piece containing the critical {\it value}. 


We call $\overline \Omega^\n$ the {\it parapuzzle piece of depth $\n\q$},
containing $\MM_\LL$.  The closure of the Misiurewicz wake $\overline \OO^\n_k$ will be also
called a parapuzzle piece of depth $\n\q+1$.
We will give now a construction of the ``parapuzzle pieces of
depth $n \geq \n\q+2$'' so that they will be the closures of the parameter domains over which     
the puzzle pieces of depth $n$ move holomorphically (with the same ``combinatorics'').
Moreover, they will form a tiling of the Misiurewicz wake, appropriately
truncated.


\begin{lem}\label{motion over O}
   Fix some Misiurewicz wake $\OO=\OO^\n_k\subset \WW$.
   Then all the boundaries of puzzle pieces up to depth $\q\n$ move holomorphically over
$\Omega^\n$, while 
   the boundaries of puzzle pieces of depth $\q\n+1$ move holomorphically over $\OO$. All these motions provide us with tubings 
   over $\OO$ respecting B\"ottcher coordinate. 
   The critical value $c\mapsto c= f_c(0)$ is a diagonal of the tubing of
   $\partial Y^{\q\n}_v$ over $\OO^\n$. Moreover, for any $c\in \OO$, 
\be\label{initial bound}
     \mod(Y^{\n\q}(c)\sm Y^{\n\q+1}(c)) \geq \de(\OO) >0.
\ee
\end{lem}

\begin{pf}
Since for $c\in \Omega^\n$,
the critical orbit $f_c^k(0)$, $k=0,1,\dots, \q\n$ does not cross the rays $\RR_i^\dyn (c)$ 
and the equipotential $\EE^\dyn_\xi(c)$, 
the configuration of Yoccoz puzzle pieces up to depth $ \q\n $ moves holomorphically
over $\Om^\n\Supset \OO$.

Similarly, for $c\in \OO$, the critical value $f_c^{\q\n+1}(0)$ does not cross 
the rays $\RR_i^\dyn (c)$  and the equipotential $\EE^\dyn_\xi(c)$ either,
so that the puzzle pieces of depth $\q\n+1$ move holomorphically over $\OO$.
By Lemma \ref{lifting of motions}, this motion is continuous up to the boundary of $\OO$.

We see that the boundary of each puzzle piece up to depth $\q\n+1$ provides us with a tubing over $\OO$. 
This tubing respects the B\"ottcher coordinate as it is induced by it.

Let us consider the puzzle piece $Y^{\q\n}_v(c)$ moving holomorphically over $\Om^\n\supset \OO$
under the B\"ottcher motion $h_c$. 
It is bounded by two arcs of external rays with some angles $\theta_+$ and $\theta_-$ 
(landing at the same point $a=a(c)$ such that $f^{\q\n}a=\alpha$), and an arc of the equipotential
$\EE_{\xi/d^{\q\n}}^\dyn$.    
By the basic relation (\ref{basic}),
the Misiurewicz wake $\OO$ is bounded by two arcs of external rays $\RR^\p_{\theta_\pm}$
(landing at the Misiurewicz root $c_*$ such that $f_{c_*}^{\q\n}(c_*)=\alpha$)
and an arc of the equipotential $\EE^\p_{\xi/d^{\q\n}}$.
Moreover, the parameter-phase map 
$$
  \gamma: \di \OO \ra \di Y^{\q\n}_v, \quad c\mapsto h_c^{-1} (c)
$$
 carries
a parameter $c\in \di \OO\sm\{c_*\} $ to the dynamical point $\gamma(c)\in  Y^{\q\n}_v\sm \{a\}$ 
with the same B\"ottcher coordinates.
This shows that the map $c\mapsto c$  satisfies properties (D1), (D2) and (D5) of the  diagonal  
to the tubing of $\partial Y^{\q\n}_v$ over $\OO$. 
It satisfies (D3) and (D4) by the discussion in \S \ref {transverse}.

Since the holomorphic motion $h$ over $ \OO $  extends to $\Om^\n$,
\be\label{delta}
   \mod (Y^0(c) \sm  Y^1_j(c) )\geq \de_\n >0, \quad  c\in \OO,
\ee
for any non-critical puzzle piece $Y^1_j$ contained in $Y^0$. 
Since 
$$
    f^{\q\n}: Y^{\q\n}(c)\sm Y^{\q\n+1}(c) \ra Y^0(c)\sm Y^1_j(c), \quad c\in \OO 
$$
is a covering of degree $d$ (for an appropriate non-critical puzzle piece $Y^1_j$),
we obtain (\ref{initial bound}). 
\end{pf}


\subsection{Puzzle motion over the parapuzzle}

Let 
$$
   \Gamma^n=\Gamma^n(c)=f^{-n}(\Gamma^0)\sm f^{-(n-1)}(\Gamma^0).
$$ 
Thus, $\Gamma^n$ is  the ``new'' ray boundary of the puzzle pieces of depth $n$ 
(which is not contained in the ray boundary of the puzzle pieces of depth $n-1$).   

We say that $f$ {\it has well defined combinatorics up to depth} $n$ if
$0$ belongs to the interior of a puzzle piece of depth $n$.  
Note that for $c\in \WW(\xi/d^\q)$ with combinatorial rotation number  ${\bf p}/\q$, 
combinatorics is well defined up to depth $\q+1$.
  
There are three ways the combinatorics can fail to be well defined at level $n$:

\ssk\nin$\bullet$
  The Julia set $J(f)$ is connected and $f^n (0)=\alpha$.
Such maps will be called $\alpha$-Misiurewicz;

\ssk\nin$\bullet$
  The Julia set is disconnected and the critical value $f^n(0)$ has height $\geq d\xi$.
This situation will be essentially avoided by appropriate shrinking of the parameter 
domains. 

 \ssk\nin$\bullet$
  The Julia set $J(f_c)$ is disconnected and $f^n(0)\in \Gamma^0$.
In this case, there are $d$ rays in $\Gamma^n$ that land at a precritical point. 
We call such precritical points {\it pinching points} 
{\it of depth} $n$ (if $n$ is the minimal integer with this property).
Note that the pinching points of depth $n$ belong to the interior of the puzzle piece
of depth $n-1$ (since by definition, they do not belong to $f^{-(n-1)}(\Gamma^0)$).

\ssk
The {\it combinatorics of $f$ up to depth} $n$
(provided it is well defined)
is the label of the puzzle piece of depth $n-1$ containing the critical value.

As we saw in \S \ref{Mis}, 
all the maps $f_c$, $c\in \OO$, have well defined combinatorics up to depth $\q\n$
(and moreover, $\OO$ is the maximal domain on which this is the case). 
We will now tile
$\overline \OO$ (truncated by appropriate equipotentials)
 according to the deeper combinatorics of the puzzle.

\begin{lem}\label{para construction}
   The set of parameters $c\in \OO$
with the same combinatorics $v$ up to depth $n \geq \q\n+1$
is an open Jordan disk bounded by the rays and equipotentials with the
same angles and heights as the puzzle piece $Y^{n-1}_v$ containing the critical value.
The closure of this disk, $\YY^n_v$, is called the parapuzzle of depth $n$ with combinatorics $v$.  
The boundaries of puzzle pieces of depth $n$
provide us with B\"ottcher tubings over
$\YY^n_v$ that fit to the tubings of the boundaries of puzzle pieces of depth $<n$ containing
it.  The critical value $\psi: c\mapsto c= f_c(0)$, $c\in \YY^n_v$, 
is a diagonal to the tubing of $\partial Y^{n-1}_v$.  
The parapuzzle pieces of depth $n$ tile the Misiurewicz wake $\OO$ truncated by 
the equipotential of height $\xi/ d^{n-1}$. 
\end{lem}

\begin{pf}
 Assume inductively that the statement is true up to depth $n$
(where the base of induction is provided by the closed Misiurewicz wake
$\overline \OO=\YY^{\q\n+1}$,
see Lemma \ref{motion over O}).
Consider one of the parapuzzle pieces, $\YY=\YY^n_j\subset \overline\OO$, and let us show how 
to tile its truncation by parapuzzle pieces of depth $n+1$. 

Let us consider the boundary of a puzzle piece $X=Y^n_k$ contained in $Y^{n-1}_v$.
By the induction assumption, it provides us with a tubing  over $\YY$
that fits to the tubing of $\partial Y^{n-1}_v$, and  $\psi$ is a diagonal to the
latter.
Hence we  can apply Lemma \ref{Arg Prin} and obtain an embedding $ \chi: \di X \ra \YY$.
The closed  disk bounded by this Jordan curve is our parapuzzle  $\XX$ of depth $n+1$.
Moreover, the map $\psi$ is the diagonal of the tubing of $\partial X$ over $\XX$.
Properties (D1), (D2) and (D5) of the diagonal follow directly from the construction,
while properties (D3) and (D4) follow from the discussion of \S \ref
{transverse}.

Since $\di X$ moves under the B\"ottcher motion and the diagonal $\psi$ is the
identity $c\mapsto c$,
the phase-parameter map $\chi: X\ra\XX$ respects the B\"ottcher coordinates.
Hence the external angles and the heights of the rays and equipotentials forming $\di \XX$ 
are the same as those of $\di X$. 

Let us now consider the puzzle pieces $Y^{n+1}_k$ of depth $n+1$.
Since $f_c(0)\not\in  \di Y^n_j(c) $ for  $c\in \XX$ and any $j$, 
these puzzle pieces move holomorphically over $\inter \XX$
(obviously, respecting the B\"ottcher coordinates).
This motion is continuous up to the boundary  by Lemma \ref{lifting of motions}.
Let us show that it fits to the tubing of the boundary of the puzzle piece of depth $n$
(and then inductively,  of all smaller depth) containing it.
Indeed, let $Y^{n+1}_l \subset Y^n_k$, and let
$\zeta\in \di Y^{n+1}_l \sm \di Y^n_k $. 
Then $\zeta\in \Gamma^{n+1}$, 
and since $\di Y^{n+1}_l$ moves under the  B\"ottcher motion,
$h_c(\zeta)\in \Gamma^{n+1}(c)$ for any $c\in \XX$.
Hence $h_c(\zeta)\not\in f_c^{-n}(\Gamma^0(c))$, 
while the latter set contains $\di Y^n_k(c)$.
This provides us with the desired tubings of depth $n+1$ over $\XX$.

Finally, the puzzle pieces $Y^n_k$ tile the puzzle piece $Y^{n-1}_v$ truncated be the
equipotential of height $\xi/d^n$ and their simultaneous motion over $\YY$ fits
to the tubing of $\di Y^{n-1}_v$. Applying Lemma  \ref{Arg Prin} once again,
we conclude that the corresponding parapuzzle pieces $\XX=\YY^{n+1}_k$ tile the puzzle piece $\YY$
truncated by the equipotential of height $\xi/d^n$.
\end{pf}

The parapuzzle piece of depth $n$ containing a point $c$ in its interior will also be denoted $\YY^n(c)$, 
(we will also use notation $\YY^n$ when the choice of the base point $c$ is self-evident or non-essential).
For instance, for $c\in \OO=\OO^\n_k$,
 we have:  
$$
   \YY^{\n\q} =\overline {\Om^\n}, \quad \YY^{1+\n\q}  = \overline {\OO^\n_k}.
$$  


For $c \in \inter \YY^n$, let 
$
  L_n= L_{n,c}: D_n \ra Y^n
$
be the {\it first landing map} to the critical puzzle piece $Y^n$.
For $n\geq \n\q+1$, its domain $D_n=D_n(c)$ consists of disjoint puzzle pieces $W_i^n=W_i^n(c)$
each of which is univalently mapped by $L_n$  onto $Y^n$.  
Note that $\C\sm \inter D_n$ is forward invariant set, and that
$\inter D_n$ contains a dense subset of the filled-in
Julia set $K(f)$.

\begin{lem}\label{Dom}
For $n \geq \n\q+1$, the set $\C\sm \inter D_n$ moves holomorphically over
$\inter \YY^n$.  This motion is equivariant, i.e., $h_c \circ f_{c_0}=f_c
\circ h_c$, and respects the B\"ottcher coordinate.
\end{lem}

\begin{pf}
  We fix some base parameter $c_0\in \YY^n$, and let $f=f_{c_0}$, $Y^n=Y^n(c_0)$, etc.  

  Let us first show  that the boundary of each domain $W=W^n_i$ moves holomorphically over $\inter \YY^n$.
Let $L_n|\, W= f^l$. Then $W$ has an itineray  $(\mu_0,\dots, \mu_{l-1},\mu_l=0)$,
where $\mu_m\not=0$ for $m<l$, 
satisfying the property: 
$$
       f^m (W)\subset Y^n_{\mu_m}, \quad m=0,\dots l.
$$
For $c=c_0$,  the restrictions $f|\, Y^n_{\mu_m}$ are univalent  
and 
$$
   f(Y^n_{\mu_m})\supset Y^n_{\mu_{m+1}}, \quad m=0,\dots, l-1.
$$ 
Since the puzzle pieces $Y^n_\mu$ move holomorphically  over $\inter \YY^n$,
the same property is valid for all $c\in \inter \YY^n$.
Now, the repeated application of Lemma \ref {lifting of motions} yields that
the boundaries of all  $f^m(W)$, $m=l-1, \dots, 0$, move holomorphically over $\inter \YY^n$
as well.

By Lemma \ref {para construction}, for
$c \in \inter \YY^n$, $f_c(0) \in \inter Y^{n-1}_v(c)$.  Thus
$0 \in \inter Y^n(c) \subset \inter D_n(c)$, and this implies that
all pre-critical points are contained in $\inter D_n(c)$.
Hence the B\"ottcher motion is well-defined on $\C \sm (\inter D_n \cup
K(f_c))$ (which is a dense subset of $\C\sm \inter D_n$).

By the $\la$-lemma, this B\"ottcher motion  
extends to  the whole set $\C\sm \inter D_n$, and this extension 
matches with the previously constructed motion of
$\displaystyle{\bigcup_i \di W^n_i}$.
The conclusion follows.
\end{pf} 

We let $h^{(n)}_c$ be the motion of $\C\sm \inter D_n$ over $\YY^n$ described in Lemma \ref{Dom}.


We say that a puzzle piece $Y=Y^k_\mu(c_0)$ {\it persists over depth} $n$
if  the boundary $\di Y^k_\mu(c)$ 
moves holomorphically over $\inter \YY^n$ respecting the B\"ottcher coordinate. 
By Lemma \ref{lifting of motions}, if $n \geq 1+\q\n$,
this motion is continuous up to the boundary of $\YY^n$
and hence provides us with a tubing over $\YY^n=\YY^n(c_0)$.
For instance, any puzzle piece $Y^n_\mu$ persists over depth $n$.

We say that a puzzle piece $Y^k_\mu$ is {\it subordinate to depth $n$} if it is not 
properly contained in some domain $W_i^n$. 

\begin{lem}\label{persistent pieces}
  Let $n\geq \n\q+1$. If a puzzle piece $Y=Y^k_\mu$
is subordinate to depth $n$, 
then  it persists over depth $n$. Moreover, its motion fits to the tubing of 
the boundary of any bigger puzzle piece $Z = Y^l_\nu$ over $\YY^n$.  
\end{lem}

\begin{pf}
The first assertion follows from Lemma \ref{Dom} 
since the boundary of $Y^k_\mu$ is contained in $\C\sm \inter D_n$.

Let us verify the second assertion.
Since ``fitting'' is a transitive property,
it is sufficient to check it for two consecutive depths, $l=k-1$.  We may
assume that $k>n$, since for $k \leq n$ the result follows from Lemma \ref
{para construction}.
Let us consider a puzzle piece $Y' = f^{k-n}(Y)$ of depth $n$,
%
and let $Z'= f^{k-n}(Z)$. The latter is a puzzle piece of depth $n-1$ containing $Y'$.

Let $h_c$ be the motion of $\C\sm \inter D_n$ from Lemma \ref{Dom}.
Since it is equivariant up to the boundary of $\YY^n$, 
we have 
\be\label{properness}
     f_c^{k-n} (h_c(\di Z)) =  h_c(\di Z'), \quad c\in \YY^n. 
\ee
  
By Lemma \ref{para construction},
$\di Y'$ provides us with a tubing over $\YY^n$ that fits to the tubing of $\partial Z'$.
By (\ref{properness}),  this property is lifted to yield that
the tubing of $\partial Y$ fits to the tubing of $\partial Z$. 
[Indeed, if $z\in \di Y\cap \inter Z$ then $f_{c_0}(z)\in \di Y' \cap \inter Z'$.
Since the tubing of $\partial Y'$ fits to the tubing of $\partial Z'$, 
$h_c(f_{c_0}^{k-n}(z))=f_c^{k-n} (h_c z)$ does not belong to $h_c(\di Z')$ for $c\in \YY^n$. 
By (\ref{properness}), $h_c z$ does not belong to $h_c (\di Z)$.]   
\end{pf}

A critical puzzle piece $Y^n$ is called a {\it child} of a critical puzzle
piece $Y^m$ ($m<n$) if the map $f^{n-m}:Y^n \to Y^m$ is unicritical.

\begin{cor} \label {nq+1}
Assume that for some $n \geq \n\q+1$  and  $k\geq 1$, the  map $f^k|\, Y^n$ is unicritical
(e.g., $Y^n$ is a child of some puzzle piece $Y^m$ and $k\in [1,n-m]$).  
Then the motion   $h^{(n)}$ provides us with a tubing of $\partial f^k (Y^n)$ over $\YY^n$, 
and the critical value $c\mapsto f_c^k (0)$ 
is a diagonal to this tubing.
\end{cor}

\begin{pf}
  The first assertion follows from Lemma \ref{persistent pieces}
since the piece $f^k (Y^n)$ is subordinate to depth $n$.
The second assertion follows from Lemma \ref {para construction} for $k=1$. 
Applying the family of univalent maps $f_c^{k-1}: f_c(Y^n(c))\ra  f^k (Y^n(c))$, 
we obtain it for any $k$. 
\end{pf}

If $f$ and $\tilde f$ have the same combinatorics up to depth $n$, a
{\it (B\"ottcher marked)
 pseudo-conjugacy} (up to depth $n$) between $f$ and $\tilde f$
is an orientation preserving homeomorphism $H:(\C,0) \to (\C,0)$
such that $H \circ f=\tilde f \circ H$ everywhere outside $\inter Y^n$,
and which is the identity near
infinity with respect to the B\"ottcher coordinates.

\begin{rem} \label {pseudocriterion}

By  Lemmas 4.2 and 4.3 of \cite {AKLS}, if $c$ and $\tilde c$ have the same
combinatorics up to depth $n$, and there exists a $K$-qc homeomorphism
$(\inter Y^n(c),0) \to (\inter Y^n(\tilde c),0)$ which is the identity at
on boundary with respect to the B\"ottcher coordinates, then
$f_c$ and $f_{\tilde c}$ are $K$-qc pseudo-conjugate (up to depth $n$).

\end{rem}

\subsection{Combinatorics of children}

If $f$ does not have well defined combinatorics of all depths, then either
the Julia set of $f$ is disconnected or the critical point is eventually
mapped to the repelling fixed point $\alpha$.  
Otherwise, we have critical puzzle pieces of all depths.
In this case, we say that $f$ is {\it combinatorially recurrent}
if the critical point returns to all critical puzzle pieces.

Given a critical puzzle piece $Y^n$, let $R_{Y^n}$ be the first return map
to $Y^n$.  The components of the domain of $R_{Y^n}$ are puzzle pieces,
which are mapped by $R_{Y^n}$ onto $Y^n$, either univalently (if the
component is non-critical), or $d$-to-$1$ (if the component is critical).
Let $\mo(Y^n)$ be the infimum, over all components $D$ of the
domain of $R_{Y^n}$, of $\mod(Y^n \setminus D)$.

If $f$ is combinatorially recurrent, then every critical puzzle piece has a child. 
 These kids are ordered by ``age'':  a child $Y^k$ is 
``older'' than a child $Y^l$ if $Y^k\supset Y^l$ (and thus, $k \leq l$).
Note that the {\it first child} $Y^k$  of $Y^n$ 
coincides with the critical component of the domain of $R_{Y^n}$.

A combinatorially recurrent map is said to be
{\it primitive renormalizable}
if there exists a critical puzzle piece $Y^n$ such that the critical point
never escapes its first child $Y^k$ under iterates of
$R_{Y^n}$: $\{R_{Y^n}^j(0) :\ j \geq 1\} \subset Y^k$.
In general, we will say that a map $f$ is {\it non-renormalizable} if it is
neither satellite nor primitively renormalizable.

A child $Q=Y^q$ of $V=Y^v$ is called {\it good} if
$f^{q-v}(0)$ is contained in the first child $U=Y^u$ of $V$.
In this situation, $K=Y^{q-v+u}$ is a child of $U$ called
the {\it friend}  of $Q$. 
Note that $f^{q-v}: Q\sm K \ra V\sm U$ is a covering map of degree $d$.  

The {\it favorite child} of $V$  is the oldest good child $Q$ that appear
after the first child $U$.  One can see that
the depth of the favorite child is the smallest
integer $q>v$ such that $f^{q-v}(0)$ belongs to the first child $U$ and
the orbit $\{f^i(0)\}_{i=1}^{q-v}$ passes through the annulus $V \setminus U$
(see the discussion preceding Lemma 2.3 of \cite {AKLS}).
If $f$ is combinatorially recurrent and non-renormalizable,
then every critical puzzle piece has a favorite child.

\subsection{Phase-parameter transfer}

 We will now apply Lemmas \ref{Arg Prin} and \ref{1}
to two dynamical situations that will often appear in what follows.

\begin{lemma} \label {12}

Let us consider four levels $s<t \leq u<w$, where $u \geq \n\q+1$,  
such that $Y^u$ is a child of $Y^s$ and $f^{u-s}(Y^w)$ is
contained in some connected component $W=W^t_i$
of the first landing map to $Y^t$.
Assume that 
\begin{enumerate}

\item $\mod(Y^s(c) \setminus W( c))>\delta$ for any $c\in \YY^u$;

\item $f_c$ and $f_{\tilde c}$ are $K$-qc pseudo-conjugate up to depth $t$ for any $c,\tl c\in Y^u$.

\end{enumerate}

Then $\mod(\YY^u \setminus \YY^w)>\eps (K,\delta)$.

\end{lemma}

\begin{pf}

Notice that $Y^s$ and $W$ are subordinate to depth $u$, 
and $c \mapsto f^{u-s}_c(0)$ is a diagonal to the tubing of $\partial Y^s$ over $\YY^u$.  
The first assumption implies that $\diam_{Y^s(c)} W(c)<M=M(\delta)$
and the second assumption implies that the maps
$h^{(t)}_c|(\partial Y^s \cup \partial W )$ have $K$-qc extensions to the whole complex plane,
$c\in \YY^u$. 
 Application of the second statement of Lemma \ref {1} gives the result:
$$
       \mod (\YY^u\sm \YY^w)\geq  \mod (\YY^u\sm \YY^{\mathrm{depth}(W)+u-s}) > \eps (K,\de).
$$
\end{pf}

\begin{lemma} \label {munu}

Let us consider four puzzle pieces $K\subset Q\subset U\ssubset V$ of respective depth
$k > q \geq u > v \geq \n\q+1$. 
Assume that $U$ is the first child of  $V$,  $Q$ is a  good child of $V$,
and $K$ is his friend. Let $\KK\subset\QQ\subset\UU\ssubset\VV$ be the corresponding
parapuzzle pieces.  
Then
$$
\mod(\QQ \setminus \KK)>\rho(\mod(\VV \setminus \QQ)) \mod(V\sm U),
$$
where $\rho:\R^+ \to \R^+$ is an increasing function.
\end{lemma}

\begin{pf}

Let $A=V\sm U$. Since this annulus persists over $\UU$, all the maps 
$$
        f_c: A\ra f(A),\quad  c\in \UU,
$$
are coverings of degree $d$.
 
By Corollary \ref{nq+1}, the puzzle piece $f(U)$ (and of course, $f(V)$) persists over  $\VV$, 
so that, the boundary of the annulus $f(A)$ moves  holomorphically under $h_c=h^v_c$.  
Let us extend this motion to the whole annulus $f(A)$
(using the same notation for the extension).
By Lemma \ref {lifting of motions},
this motion lifts to a holomorphic motion $H_c$ of $A$ over $\inter \UU$ continuous up to the boundary.
For any $z \in A \setminus \partial V$ and $c \in \UU$,
$h_c(f(z)) \notin \partial f(V)$, hence $H_c(z) \notin \partial V$.
Thus, the motion of $A$  fits to the tubing of $\di V$ over $\UU$.

Let us now consider a unicritical family 
$f_c^{q-v} : Q(c) \ra V(c)$ over $\QQ$ and a unicritical family
$f_c^{q-v} : K(c) \ra U(c)$ over $\KK$.
By Corollary~\ref{nq+1}, the critical value $c \mapsto f^{q-v}_c (0)$ is a diagonal to the corresponding tubings:
of $\partial V$ over $\QQ$ and of $\partial U$ over $\KK$. 
Hence the corresponding  phase-parameter map $\chi: V\ra \QQ$ maps the annulus $A$
onto the annulus $\QQ \sm \KK$. 
By Lemma~\ref {Arg Prin}, the dilatation of this map is bounded by the dilatation of the motion 
$H_c$ over $\QQ$, which is equal to the dilatation of $h_c|\, f(A)$ over $\QQ$. By the $\la$-lemma, 
the latter is bounded by $\rho(\mod(\VV \setminus \QQ))$,
which implies the desired estimate.
\end{pf}

\section{The favorite nest and the principal nest}

\subsection{The favorite nest} \label {Kn}

Let $Q^0=Y^{\n\q}$.  Let $Q^{i+1}$ be the favorite child of $Q^i$,
and let $P^i$ be the first child of $Q^i$.  
Let $q_n$ and $p_n$ be the depths of these puzzle pieces, i.e., 
$Q^n=Y^{q_n}$, $P^n=Y^{p_n}$.  
Let $k_n=q_n+p_{n-1}-q_{n-1}$  and $K^n=Y^{k_n}$.
Note that 
$$
   f^{q_n-q_{n-1}} : (Q^n, K^n) \ra (Q^{n-1}, P^{n-1}).
$$
By Proposition 2.4 of \cite {AKLS}, we have
\be \label {Qn-1P}
\mod(Q^n \setminus P^n)>\delta(\OO), \quad n \geq 0.
\ee
This implies in particular that
\be \label {Qn-1K}
\mod(Q^n \setminus K^n)=\frac {1} {d} \mod(Q^{n-1} \setminus P^{n-1})>\frac
{1} {d} \delta(\OO), \quad n \geq 1.
\ee

Let us also consider the corresponding parapuzzle pieces:
$\QQ^n=\YY^{q_n}$, $\PP^n=\YY^{p_n}$, and $\KK^n=\YY^{k_n}$.

\begin{thm} \label {favorbounds}

There exists $\delta>0$, depending only on $\OO$, such that
for $n \geq 2$, $\mod(K^n \setminus P^n)>\delta$,
$\mod(\KK^n \setminus \PP^n)>\delta$,
and for $n \geq 3$, $\mod(\QQ^n \setminus \KK^n)>\de$.

\end{thm} 

\begin{pf}

We start with the first and second estimates.
The map $f^{q_n-q_{n-1}}$ is unicritical on $Q^n$; 
all the more, it is unicritical on $P^n$. It follows that 
$p_n-q_n \geq q_n-q_{n-1}$, and hence the puzzle piece
$$
    D = f^{q_n-q_{n-1}}(P^n)\subset P^{n-1}
$$
is a component of the first landing map to $Q^n$.
It follows that $f^{p_{n-1}-q_{n-1}}(D)$ is contained in a component of the
first return map to $Q^{n-1}$.  In particular
\begin{align} \label {Pn-1D}
\mod(K^n \setminus P^n)&=\frac {1} {d}
\mod(P^{n-1} \setminus D) \geq \frac {1} {d^2} \mo(Q^{n-1})\\
\nonumber
&\geq \frac{1}{d^3}\mod(Q^{n-2}\sm P^{n-2})
    \geq \frac{1}{d^3} \de(\OO),
\end{align}
where the last estimate follows from (\ref {Qn-1P}), 
while the previous one follows from Lemma 2.2 of \cite {AKLS}.
This proves the first estimate of the lemma. 

For the second estimate, let us define $s<t<u<w$ as follows: $s=p_{n-1}$,
$t=q_n$, $u=k_n$, $w=p_n$.  Then (\ref {Pn-1D}) 
implies that condition (1) of Lemma \ref {12}
is satisfied. 
Condition (2) of Lemma \ref {12}
is satisfied 
by Theorem 4.4 (and Remark~4.1)  of \cite {AKLS}.  
Applying Lemma \ref {12}, we get the conclusion.

\msk
For the third estimate, let $s<t<u<w$ be as follows:  $s=q_{n-1}$,
$t=k_{n-1}$, $u=q_n$, $w=k_n$.  Then
$f^{u-s}(Y^w)=f^{q_n-q_{n-1}}(K^n)=P^{n-1}$
is contained in $K^{n-1}$,
which is a (trivial)
component of the first landing map to $K^{n-1}$.  By (\ref
{Qn-1K}),
condition (1) of Lemma \ref {12} is satisfied.

Furthermore, $\mod(\KK^{n-1}
\setminus \QQ^n) \geq \mod(\KK^{n-1} \setminus \PP^{n-1}) \geq \delta$ by
the previous estimate.  By Remark \ref {pseudocriterion} and
the $\la$-lemma,
this implies that condition (2) of Lemma~\ref {12} is
satisfied with $K=K(\OO)$.  The conclusion now follows from Lemma~\ref {12}.
\end{pf}

\subsection{The principal nest}

Let $V^0=Y^{\n\q}$.
The {\it principal nest} starting at $V^0$ is the nest
$V^0 \supset V^1 \supset ...$ obtained by taking
$V^{i+1}$ as the first child of $V^i$.  
Let $v_i$ be the depth of $V^i$, i.e.,  $V^i=Y^{v_i}$.
Let $\VV^i=\YY^{v_i}$ stand for the corresponding parapuzzle pieces.
By Corollary \ref{nq+1}, 
$\partial V^i$ provides us with a tubing over $\VV^{i+1}$ with diagonal  $c\mapsto f_c^{v_{i+1} - v_i }(0)$

Let $g_i=g_{i,c}=R_{V^i}$. 
If the critical point returns to $V^{i+1}$, 
we let $s_i=s_i(c)$ be the smallest $k \geq 0$ such that $g_i^{k+1}(0) \in V^{i+1}$.
(In other words, $s_i+1$ is the first return
time of the critical point back to $V^{i+1}$ under the iterates of $g_i$.)  
If $s_i=0$ 
(that is, $g_i (0) \in V^{i+1})$, we say that the return to $V^i$ is {\it central}.

The map $g_{i+1}=f^{v_{i+2}-v_{i+1}}: V^{i+2}\ra V^{i+1}$ admits a unicritical  ``Koebe extension''
$$
   \Bg_{i+1}=\Bg_{i,c}: \BV^{i+2}\ra V^i,\quad 
   {\mathrm {where}} \quad \BV^{i+2}= Y^{v_{i+2}-v_{i+1}+v_i}\equiv Y^{\Bv_{i+2}}\subset V^{i+1}.
$$
(Note that $\BV^{i+2}$ is a good child of $V^i$.)
By Corollary \ref{nq+1}, $\partial V^i$ provides us with a tubing over
the parapuzzle $\tVV^{i+2}\equiv Y^{\Bv_{i+2}}$
with diagonal $c\mapsto f_c^{v_{i+2}-v_{i+1}}(0)$.

\begin{lemma}\label{6.1}

We have the estimates
\be \label {6.11}
\mod(V^{i+1} \setminus V^{i+2}) \geq \frac {1} {d} (\mod(V^i \setminus
V^{i+1})+s_i \mo(V^i)),
\ee
\be \label {6.12}
\mod(\VV^{i+1} \setminus \VV^{i+2}) \geq \mod(\tVV^{i+2} \setminus \VV^{i+2})> 
    \rho(\mod(\VV^i \setminus  \VV^{i+1})) \mod(V^i \setminus V^{i+1}).
\ee

\end{lemma}

\begin{pf}

Let $D^j=Y^{d_j}$, $1 \leq j \leq s_i+1$ be the pullback of $V^i$ under
$g_i^j$, so that
\be\label{D's}
      V^{i+1}=D^1 \supset ... \supset D^{s_i+1} = \BV^{i+2}. 
\ee  
Notice that the $D^j$ are all children of $V^i$.

For $1 \leq j \leq s_i$, $f^{d_j-v_i}(D^{j+1})$
is a non-central component of the $\Dom g_i$.  Thus
$\mod(V^i \setminus f^{d_j-v_i}(D^{j+1})) \geq \mo(V^i)$.
Since 
$$
   f^{d_j-v_i}:D^j \setminus D^{j+1} \to V^i \setminus f^{d_j-v_i}(D^{j+1})
$$
is a covering of degree $d$, we have
\be \label {6.14}
\mod(D^j \setminus D^{j+1}) \geq \frac {1} {d} \mo(V^i), \quad
1 \leq j \leq s_i.
\ee
Moreover,
$f^{v_{i+2}-v_{i+1}}: \BV^{i+2}\sm V^{i+2} \to V^i\sm V^{i+1}$ is a covering of degree $d$ as well.  
Hence
\be \label {6.15}
\mod(\BV^{i+2} \setminus V^{i+2})=\frac {1} {d} \mod(V^i \setminus V^{i+1}).
\ee
Putting  (\ref {6.14}), (\ref {6.15}) together with the Gr\"otzcsh inequality, we get (\ref {6.11}).

Applying to the nest  $V^{i+2}\subset \BV^{i+2}\subset V^{i+1} \subset V^i$  Lemma \ref {munu}, 
we obtain
$$
\mod(\tVV^{i+2} \setminus \VV^{i+2})>\rho(\mod(\VV^i \setminus \tVV^{i+2})) \mod(V^i \setminus V^{i+1}).
$$
Since $\VV^{i+1} \supset \tVV^{i+2}$,  (\ref {6.12}) follows.
\end{pf}

Define $0=i_0<i_1<...$ as the sequence of levels such that for $j>0$ the return to $i_j-1$ is non-central,
i.e., $g_{i_j-1}(0)\not\in V^{i_j}$.

\begin{lemma}\label{6.2}

Let $i_{j-1} \leq s<i_j$.  Then we have the estimates
$$
\mod(V^{i_j} \setminus V^{i_j+1}) \geq
\frac {1} {d} \mod(V^s \setminus V^{s+1}),
$$
$$
\mod(\VV^{i_j} \setminus \VV^{i_j+1})>\rho(\mod(\VV^s \setminus \VV^{s+1}))
\mod(V^s \setminus V^{s+1}).
$$

\end{lemma}

\begin{pf}

Let $n$ be the first moment such that $f^n(0) \in V^s
\setminus V^{s+1}$, and let $m>0$ be the first moment such that
$f^{n+m}(0) \in V^{s+1}$.  Then
$Q=Y^{v_s+n+m}$ is the favorite child of $V^s$. 

We have, $Q\supset V^{i_j+1}$
since $V^s\supset V^{i_j}$ and $n+m$ is not bigger than the first return time to $V^{i_j}$.
A similar argument gives $f^{n+m}(V^{i_j+1}) \subset V^{s+1}$.
Hence
$$
   \mod(Q \setminus V^{i_j+1}) \geq \frac {1} {d} \mod(V^s \setminus V^{s+1}).
$$  
Since  $Q\subset Y^{v_s+n}=V^{i_j}$,
the first statement follows.

The second statement follows from Lemma~\ref {munu} applied to the nest 
$$
   K\subset Q\subset V^{s+1}\subset V^s,
$$
where $K$ is the friend of $Q$ (since $f^{n+m}(V^{i_j+1}) \subset V^{s+1}$
implies $V^{i_j+1} \subset K$).
\end{pf}

\begin{lemma} \label {QVQ}

If $V^{i_j-1} \subset Q^n$ and
$V^{i_j+1}$ is defined then $V^{i_j+1}
\subset Q^{n+1}$.

\end{lemma}

\begin{pf}

Recall that $q_{n+1}=q_n+m$
where $m$ is minimal with $f^m(0) \in P^n$ and
$\{f^k(0) :\ 1 \leq k \leq m\} \cap Q^n \setminus P^n \neq \emptyset$.
Clearly $V^{i_j} \subset P^n$, so
we just have to show that $\{f^k(0) :\ 1 \leq k \leq
v_{i_j+1}-v_{i_j}\} \cap Q^n \setminus P^n \neq \emptyset$, as this implies
$q_{n+1}=q_n+m \leq v_{i_j}+m \leq v_{i_j}+v_{i_j+1}-v_{i_j}$.


Let $k \geq 1$ be minimal such that $f^k(0) \in Q^n \setminus P^n$.
Then $k=l_0(p_n-q_n)$ for some $l_0>0$, and
$R_{Q^n}^l(0) \in Y^{q_n+(l_0-l)(p_n-q_n)} \setminus
Y^{p_n+(l_0-l)(p_n-q_n)}$, $1 \leq l \leq l_0$.

Since $R_{V^{i_j-1}}(0) \notin V^{i_j}$, we also have
$R_{Q^n}(0) \notin V^{i_j}$, so that $V^{i_j} \subset
Y^{q_n+(l_0-1)(p_n-q_n)}$.  This clearly implies that
$R_{Q^n}^l(0) \notin V^{i_j}$, $1 \leq l \leq l_0$, so $v_{i_j+1}-v_{i_j}
\geq k$, as desired.
\end{pf}

\begin{thm} \label {bounds}

There exists $\delta>0$, depending only on $\OO$ such that
\begin{enumerate}
\item $\mod(V^{i_j} \setminus V^{i_j+1})>\delta$,
\item $\mod(\VV^{i_j} \setminus \VV^{i_j+1})>\delta$.
\end{enumerate}

\end{thm}

\begin{pf}

By Lemma \ref {motion over O}, the first statement holds for $j=0$.  Since
$\VV^0=\YY^{\n\q}=\Omega^\n$ and $\VV^1 \subset \YY^{\n\q+1}=\OO$, we have
$\mod(\VV^0 \setminus \VV^1) \geq \mod(\Omega^\n \setminus
\OO)=\delta(\OO)$, so the second statement also holds for $j=0$.

By Lemmas \ref {QVQ} and \ref{6.2},
it is enough to show that for every $n \geq 3$ such that $Q^{n+1}$ is
defined, if $k \geq 0$ is maximal such that $V^k \supset Q^n$ then
$\mod(V^{k+1} \setminus V^{k+2}),\mod(\VV^{k+1} \setminus
\VV^{k+2})>\delta$.

Note first that since $P^n$ and $V^{k+2}$ are the first children of 
$Q^n$ and $V^{k+1}$ respectively, we have: $V^{k+2}\subset P^n$.
Recall the definition of $K^n$ given in the beginning of section \ref {Kn}:
$K^n=Y^{q_n+p_{n-1}-q_{n-1}}$.
If $V^{k+1} \supset K^n$, then 
$$ \mod(V^{k+1} \setminus V^{k+2}) \geq \mod(K^n \setminus P^n)\quad  \mathrm{and} \quad 
  \mod(\VV^{k+1} \setminus \VV^{k+2}) \geq \mod(\KK^n \setminus \PP^n).
$$
If $K^n \supset V^{k+1}$, then 
$\mod(V^k \setminus V^{k+1}) \geq \mod(Q^n \setminus K^n)$ and
$\mod(\VV^k \setminus \VV^{k+1}) \geq \mod(\QQ^n \setminus \KK^n)$,
so by Lemma \ref {6.1},
$$
\mod(V^{k+1} \setminus V^{k+2}) \geq \frac {1} {d}
\mod(Q^n \setminus K^n),
$$
$$
\mod(\VV^{k+1} \setminus \VV^{k+2})>\rho(\mod(\QQ^n \setminus \KK^n))
\mod(Q^n \setminus K^n).
$$
In either case, the result follows from Theorem \ref {favorbounds}.
\end{pf}

\begin{rem}

In \cite {puzzle}, \cite {parapuzzle}, it is shown that if $d=2$
then one has better estimates
\begin{enumerate}
\item $\mod(V^{i_j} \setminus V^{i_j+1})>(j+1)\delta$,
\item $\mod(\VV^{i_j} \setminus \VV^{i_j+1})>(j+1)\delta$.
\end{enumerate}

\end{rem}

\begin{rem}[{\it Beau bounds}] \label {beau bounds}

It follows from the arguments in
\cite {KL1}, \cite {KL2}, \cite {AKLS}, and this
work that there exists $\delta>0$ (depending on the degree, but not on
$\OO$) such that for every
$j$ sufficiently large (depending on $\OO$) one has
\begin{enumerate}
\item $\mod(V^{i_j} \setminus V^{i_j+1})>\delta$,
\item $\mod(\VV^{i_j} \setminus \VV^{i_j+1})>\delta$.
\end{enumerate}

\end{rem}

%% file: statistics.tex
\section{Slow recurrence}

If $c \in \MM$ is not combinatorially recurrent then either
$c$ has a non-repelling fixed point, or $f_c$ is satellite renormalizable,
or $f_c$ is {\it semi-hyperbolic}
(that is, its critical point is non-recurrent and belongs to the Julia set).
It is well known that the set of semi-hyperbolic parameters has zero
Lebesgue measure.  Indeed, in \cite {RL} a more precise version of the
following is proved:

\begin{thm}\label{Misiur}

If $c \in \MM$ is a semi-hyperbolic
parameter then $c$ is a Lebesgue density point of the complement of $\MM$.

\end{thm}

In particular, almost every parameter in some $\VV^{n+1}$ is either in the
complement of $\MM$ or is combinatorially recurrent.  For real parameters,
the corresponding statement has been proved by \cite {Sa}: the set of
semi-hyperbolic parameters $c \in \MM \cap \R$ has zero one-dimensional
Lebesgue measure.  We will now
concentrate on the analysis of combinatorially recurrent parameters.

\begin{rem}

A proof that the set of semi-hyperbolic parameters
has zero Lebesgue measure can be also
obtained along the lines of the approach we follow for combinatorially
recurrent parameters given below. 

\end{rem}


\subsection{Positive probability of slow recurrence}


\comm{
Let
$$
\tVV^n=\tVV^n(c_0) = \YY^{v_n-v_{n-1}+v_{n-2}}(c_0), \quad n \geq 2.
$$
The unicritical  family $f_c^{v_n-v_{n-1}}: Y^{v_n-v_{n-1}+v_{n-2}}(c)\ra V^{n-2}(c)$ 
over $\tVV^n(c_0)$
provides us with a ``Koebe extension'' (in both dynamical and parameter directions) 
to the unicritical family  $f_c^{v_n-v_{n-1}}: V_c^n\ra V_c^{n-1}$ over $\VV^n$. 
Let $g_{n-2,c}$ be the return map to $V^{n-2}(c)$ under iterates of
to $f_c$. 
By Corollary \ref {nq+1}, the critical value
$$
  c\mapsto f_c^{v_n-v_{n-1}}(0) = g_{n-2,c}^{s_{n-2}+1}(0), 
    \quad \mathrm{where}\quad v_i=v_i(c_0),\ s_{n-2}=s_{n-2}(c_0),  
$$
is a diagonal to the tubing of $\partial V^{n-2}$ over $\tVV^{n-2}$
and is a diagonal to the tubing of $\partial V^{n-1}$ over $\VV^n$ . 
It follows that  $s_{n-2}(c)\geq s_{n-2}(c_0)$ for $c\in \tVV^n$,
with equality attained iff $c\in \VV^n(c_0)$. 
Notice also that  $s_k(c)= s_k(c_0)$ for $k<n-2$ and
$c\in \tVV^n(c_0)$.   
******************}

\begin{lemma}\label{mod growth}

For any $\de>0$ and $\la>0$
there exists $\epsilon=\epsilon(d,\delta)$ and
$K= K(d,\delta,\la)$ with the following property.
Assume that for some parapuzzle piece $\VV^n$,  
$$ 
  \mod(\VV^n \setminus \VV^{n+1})>\delta  \quad {\mathrm{and}} \quad \mod(V^n(c) \setminus V^{n+1}(c))>\delta, \quad c\in \VV^{n+1}.
$$ 
Then for every $c \in \VV^{n+k+3}$ such that $s_{n+i} \geq (1+i) K$, $i=0,1,\dots, k$,  we have:
\begin{enumerate}
\item $\mod(V^{n+i} \setminus V^{n+i+1})>\max \{\epsilon,(i-1)\lambda\}$,
\item $\mod(\VV^{n+i} \setminus \VV^{n+i+1})> \max \{\epsilon,(i-2) \la\}$, 
\item $\mod(\tVV^{n+i+2} \setminus \VV^{n+i+2})> \max \{\epsilon,(i-1) \la\}$, 
\end{enumerate}
for  $i=0,1,\dots, k+1$.
\end{lemma}

\begin{pf}

The inequalities
$\mod(V^{n+i+1} \setminus V^{n+i+2}) \geq \frac {1} {d} (1+s_{n+i})
\mo(V^{n+i})$, given by (\ref {6.11}), and
$\mo(V^{n+i+1}) \geq \frac {1} {d}
\mod(V^{n+i} \setminus V^{n+i+1})$, given by Lemma 2.2 of \cite {AKLS},
imply the estimate
$\mod(V^{n+i} \setminus V^{n+i+1})>\delta_i (1+i)$, where
$\delta_0=\delta$, $\delta_1=\delta/d$, and $\inf_{i \geq 1}
\delta_i$ goes to infinity with $K$ (given $\delta$ and $d$ fixed).  The
first estimate follows.

Together with (\ref {6.12}), it implies the rest. 
\comm{*************
By (\ref {6.12}), we also have:
$$
\mod(\VV^{n+i+1} \setminus \VV^{n+i+2})>\rho(\mod(\VV^{n+i}
\setminus \VV^{n+i+1})) \mod(V^{n+i} \setminus V^{n+i+1}),
$$
where
$\rho:\R^+ \to \R^+$ is some increasing function.  The second estimate
follows.

\msk
Let us now consider a unicritical family 
$$
   f_c^{v_{n+i+2}-v_{n+i+1}}: V^{n+i+2}(c)\ra V^{n+i+1}(c) \quad \mathrm{over}\quad  \VV^{n+i+2}
$$
 and its Koebe extension   
$$
   f_c^{v_{n+i+2}-v_{n+i+1}}: Y^{v_{n+i+2}-v_{n+i+1}+v_{n+i}}(c)\ra V^{n+i}(c) \quad \mathrm{over}\quad  \tVV^{n+i+2}.
$$
As we know,  the critical value $c\mapsto f_c^{v_{n+i+2}-v_{n+i+1}}(0)$ is a diagonal to the both families.
Moreover, the restriction of $f$ to the annulus  $A=V^{n+i}\sm V^{n+i+1}$ is a $d$-to-1 covering onto the image,
and by Lemma \ref{persistent pieces} the puzzle piece  $f(V^{n+i+1})$ is subordinate to depth $v_{n+i}$.  
By Lemma \ref {munu} (with $u=v_{n+i+2}-v_{n+i+1}+v_{n+i}$, $t=w=v_{n+i}$,
$r=1$, $Y^k_\mu=V^{n+i}$, and $Y^l_\nu=V^{n+i+1}$),  
$$
\mod(\tVV^{n+i+2} \setminus \VV^{n+i+2})>
\rho(\mod(\VV^{n+i} \setminus \tVV^{n+i+2})) \mod(V^{n+i} \setminus V^{n+i+1})
$$
$$
\geq \rho(\mod(\VV^{n+i} \setminus \VV^{n+i+1})) \mod(V^{n+i} \setminus V^{n+i+1}),
$$
and the third estimate follows. 
***********************}
\end{pf}

Given a sequence of disjoint sets $X^n$, $n \geq 0$,
on a probability space $\VV$,
we let $\X^n =\cup_{k<n} X^n$, $\X\equiv \cup_{k \geq 0} X^k$. 
Below we will make use of the following general formula: 
\be\label{P}
  1-\P(\X) = \prod_{n \geq 0} (1-\P(X^n|\,\VV\sm \X^n))
\ee
(where $\P$ stands for probability or conditional probability).
Indeed, letting $A^n= \VV\sm \X^n$, $A\equiv \cap_{n \geq 0} A^n=\VV
\setminus X$,   
we have: $A^0 \supset A^1\supset\dots \supset A$, and
$$
  \P(A) = \prod_{n \geq 0} \P(A^{n+1}|\, A^n)= \prod_{n \geq 0}
\P(\VV\sm X^n|\, \VV\sm \X^n),
$$
which yields (\ref{P})

\begin{lemma}\label{linear growth}

Assume that for some parameter $c_0$, 
$$ 
   \mod(\VV^n \setminus \VV^{n+1})>\delta, \quad \mod(V^n \setminus V^{n+1})>\delta.
$$  
Let $Z_{nr} \subset \VV^{n+1}$ be the
set of parameters which are not combinatorially recurrent.
Fix some $K>0$ as in Lemma \ref{mod growth},
and let
$Z_{sr} \subset \VV^{n+1}$ be the set of combinatorially recurrent parameters
for which $s_{n+k} \geq (1+k)K$, $k \geq 0$.%
\footnote{Label ``$sr$'' stands for ``slow recurrent''.}

  Then
$\P (Z_{sr} \cup Z_{nr}\, |\, \VV^{n+1})>\epsilon(\delta,d,K)>0$.

\end{lemma}

\begin{pf}

We can assume that $K$ is larger than the $K(\delta,d)$ given by
Lemma \ref{mod growth}.
Let $t_k= (1+k)K $, and  let $\X=\VV^{n+1} \setminus (Z_{sr} \cup Z_{nr})$.
For $k\geq 0$, $0\leq j<t_k$, let $X^{k,j} \subset \X$
be the set of all $c\in \X $ such that
$s_{n+i} \geq t_i$, $0 \leq i<k$ and $s_{n+k}=j$.  
Notice that $\displaystyle{ \X=\bigsqcup_{(k,j)} X^{k,j}}$. 
We order the pairs $(k,j)$ lexicographically. As above, let
$$
   \X^{k,j}=\bigcup_{(k',j')<(k,j)} X^{k',j'}
$$
Notice that for $c \in X^{k,j}$, we have:
\be\label{note1}
      \tVV^{n+k+2}(c) \cap \X^{k,j}=\emptyset
\ee
while
\be\label{note2} 
      \tVV^{n+k+2}(c)\cap  X^{k,j} \subset \VV^{n+k+2}(c).
\ee
Indeed, for $\tilde c \in \tVV^{n+k+2}(c)$, we have: 
$s_{n+i}(\tilde c) = s_{n+i}(c)\geq t_{n+i}$ for $i< k$, while
$s_{n+k}(\tilde c)\geq s_{n+k}(c) = j$, 
with equality attained iff $\tilde c \in \VV^{n+k+2}(c)$.

Together with (\ref{P}), (\ref{note1}) gives us:
$$
                \P(Z_{nr}\cup Z_{sr}|\, \VV^{n+1})  =  1-\P(\X|\, \VV^{n+1})
$$
$$
 = \prod_{(k,j)} ( 1- \P(X^{k,j}|\, \VV^{n+1} \setminus \X^{k,j}))
  \geq \prod_{(k,j)} (1-\sup_{c \in X^{k,j}} \P(X^{k,j}|\, \tVV^{n+k+2}(c)).
$$
It is thus enough to prove an estimate such as
$$
   \P(X^{k,j}|\tVV^{n+k+2}(c)) \leq e^{-(1+k)\eps},\quad c \in X^{k,j}
$$
for some $\epsilon=\epsilon(\delta,d)$.
But this follows from (\ref{note2}) and the
estimate 
$$\mod(\tVV^{n+k+2}(c) \setminus \VV^{n+k+2}(c)) \geq (1+k)\eps $$
 of the previous lemma.
\end{pf}

\begin{rem}

The above proof can be easily refined as follows.  
One can define $Z_{sr}$ as the set of
combinatorially recurrent parameters $ c \in \VV^{n+1}$ for which the
sequence $s_{n+i}$ satisfies:
$s_n,s_{n+1},s_{n+2} \geq K$ and $s_{n+i+1} \geq 2^{s_{n_i}}$
for $i \geq 2$ (thus displaying
``torrential growth'' in the terminology of \cite {AM1}).
We would still obtain
$\P(Z_{nr}\cup Z_{sr})|\, \VV^{n+1}) > \epsilon$.

\end{rem}

Let $\SSS \subset \MM$ be the set of combinatorially recurrent parameters
$c$ such that $s_n$, $\mod(\VV^n(c) \setminus \VV^{n+1}(c))$ and
$\mod(V^n(c) \setminus V^{n+1}(c))$ grow at least linearly with $n$.

Let $Z_r$ be the set of combinatorially recurrent non-renormalizable
parameters in $\MM$. 

\begin{cor} \label {SSS}

For $c\in Z_r$,  
there exist parapuzzle pieces $\VV^{n+1}(c)$ of
arbitrarily small diameter such that 
$\P(\MM \setminus \SSS|\VV^{n+1}(c) )<1-\delta$, with $\delta=\delta(\OO)$.%
\footnote {Actually $\delta>0$ does not depend on $c$, not even via
the Misiurewicz wake, see Remark~\ref{beau bounds}.}

\end{cor}

\begin{pf}

For $c\in Z_r$, the sequence $i_j$ in Theorem \ref {bounds} is infinite.  
By the Rigidity  Theorem of \cite{AKLS} (or directly from Theorem \ref {bounds}),
the  parapuzzle pieces $\VV^{i_j}(c)$ shrink to $c$.
We can now apply Lemma \ref{linear growth} with $n=i_j$, 
which implies the statement (since By Theorem \ref{Misiur}
combinatorially non-recurrent parameters in $\VV^{n+1}$ are almost surely outside $\MM$).
\end{pf}

In order to exploit the previous corollary, we will need the following
``Density Points Argument''.
Let us consider a measurable set $X \subset \C$ such that for almost every
$x \in X$ there exists a sequence $X^n(x) \subset \C$
of measurable sets containing $x$ such that $\diam X^n(x) \to 0$.  Assume
that any two $X^n(x)$, $X^m(y)$ are either nested or disjoint.  Then
$\lim \P(X\,|\,X^n(x))=1$ for almost every $x \in X$.  This is a particular
case of the standard generalization of the Lebesgue Density Points Theorem
(which assumes that the family $\{X^n(x)\}_{x,n}$ satisfies the
Besikovic Covering Property), and can be also seen as a
direct consequence of the Martingale Convergence Theorem.

\begin{cor} \label {aesss}

For almost every $c \in \MM$, either $f_c$ has an attracting fixed point, or
$f_c$ is renormalizable, or $c \in \SSS$.

\end{cor}

\begin{pf}

It is enough to show that $\SSS$ has full Lebesgue
measure in $Z_r$.  For fixed $n$, the parapuzzle pieces $\VV^n(\tl c)$ define a
partition of $Z_r$.  Since the $\VV^n(c)$ shrink to $c$ for any 
$c \in Z_r$, we can apply the
Density Points Argument, which implies that
$\lim \P(Z_r \setminus \SSS\, |\, \VV^n(c))=1$ for
almost every $c \in Z_r \setminus \SSS$.
But by Corollary~\ref{SSS},  this can not happen for $c \in Z_r$.
\end{pf}

\subsection{Real parameters}

Our entire discussion goes through for real parameters as well, without
changes.  However there are no parameters in $\VV^{n+1} \cap \R \setminus
\MM$, so that we can state the following stronger version of
Corollary~\ref {SSS}:

\begin{cor}

There exists $\delta>0$ such that for $c \in Z_r \cap \R$,
there exist parapuzzle pieces $\VV^{n+1}(c)$ of
arbitrarily small diameter such that $\P(\SSS|\VV^{n+1}(c) \cap \R)>\delta$.

\end{cor}

\begin{cor} \label {aesss2}

For almost every $c \in \MM \cap \R$, either $f_c$ has an attracting
fixed point, or $f_c$ is renormalizable, or $c \in \SSS$.

\end{cor}

Parameters in $\SSS \cap \R$ have {\it exponential decay of geometry}, that
is, the ratios $\lambda_n$ between the lengths
of $V^{n+1} \cap \R$ and $V^n
\cap \R$ satisfy $\lambda_n<C e^{-\epsilon n}$ for some $C>0$, $\epsilon>0$.
Hence 
\be
\sum \lambda_n^{1/d}<\infty,
\ee
and by the Martens-Nowicki Criterium \cite {MN}
the maps $P_c$, $c\in \SSS$, are {\it stochastic} 
(that is, they have an absolutely continuous invariant
measure).

\begin{rem}

In \cite {BSS1} it is shown that decay of geometry (that is,
$\lambda_n \to 0$) already implies the existence of an absolutely
continuous invariant measure.

\end{rem}

\begin{cor} \label {j}

The set of non-renormalizable
stochastic parameters $c \in \MM\cap \R$ has positive Lebesgue measure.

\end{cor}

\begin{cor} \label {p}

Almost every non-renormalizable parameter $c \in \MM \cap \R$ is either regular or
stochastic.

\end{cor}

\begin{rem}

Corollary \ref {j}, in the case $d=2$, was obtained in \cite {J}.  The
generalization to the higher degree case is well known
(see \cite {T}, Theorem 2, which follows the approach of \cite {BC}).
Our proof is rather different.

Corollary \ref {p}, in the case $d=2$, was obtained in \cite {parapuzzle}
and is new in the higher degree case.

\end{rem}

%% file: outline.tex
\section{Conclusion}

By now, we have carried out all the extra work needed for the higher degree
case: once we know that the phase-parameter geometry almost surely decays
(Corollaries \ref{aesss} and \ref {aesss2}),
the further argument is the same as in the quadratic one.   
For reader's convenience, below we will briefly elaborate this statement.

\subsection{Real parameters}

\subsubsection{Collet-Eckmann property (Theorem \ref {ce})}

Standard renormalization considerations reduce the analysis of exactly
$n$-times renormalizable parameters with some fixed combinatorics to the
analysis of non-renormalizable parameters in a ``Multibrot-like
family''.  The analysis of Multibrot-like families is parallel to the one
we have done (see \cite {parapuzzle} which deals directly, in the case
$d=2$, with Mandelbrot-like families), and one reaches the same theorems,
with the difference that all constants may depend on the geometry of the
Multibrot-like family under consideration.

Since a renormalizable map is Collet-Eckmann if and only if its
renormalization is, Theorem \ref {ce} follows from the statement that (in a
Multibrot-like family) real non-renormalizable parameters are almost
surely either regular or Collet-Eckmann.
In view of Corollary \ref{aesss2}, this is
reduced to the following result:

\begin{thm}

Collet-Eckmann parameters have full (one-dimensional) Lebesgue measure
in $\SSS \cap \R$.

\end{thm}


This result follows from the statistical argument of \cite{AM1}:
as it is pointed out in Remark 2.1 of that paper, 
the statistical argument applies in any degree case to the
set of parameters satisfying the following properies:

\ssk\nin $\bullet$ $\liminf s_n \geq 1$;

\ssk\nin $\bullet$ Exponential decay of the real phase geometry 
(meaning that the ratios of the lengths of the real traces of  $V^{n+1}$ and $V^n$
decay exponentially);

\ssk\nin $\bullet$
   Growth of the parameter moduli $\mod(\VV^n \setminus \VV^{n+1})$. 

\ssk
All these conditions hold for non-renormalizable parameters $c \in \SSS$ 
(the exponential decay of geometry
follows from the linear growth of the phase moduli $\mod(V^n \setminus V^{n+1})$).

\subsubsection{Further statistical properties}

The statistical analysis of \cite {AM1}
and \cite {AM4} goes far beyond the Collet-Eckmann property, and
gives a very detailed description of
maps in $\SSS \cap \R$.
As for the Collet-Eckmann property, it can be directly applied to the higher degree
case:

\begin{thm}

For almost every $c \in \R$ such that $f_c$ is not regular or infinitely renormalizable,
\begin{enumerate}
\item The critical point is polynomially recurrent with exponent $1$:
$$
\limsup \frac {-\ln |f_c^n(0)|} {\ln n}=1,
$$
\item The critical orbit is equidistributed with respect to the absolutely
continuous invariant measure $\mu$:
$$
\lim \frac {1} {n} \sum_{i=0}^{n-1}
\phi(f_\lambda^i(0))=\int \phi d\mu
$$
for any continuous function $\phi:I \to \R$,
\item The Lyapunov exponent of the critical value,
$\displaystyle{\lim \frac {1} {n} \ln |Df^n(f(0))|}$,
exists and coincides with the Lyapunov exponent of $\mu$.
\item The Lyapunov exponent of any periodic point $p$ contained in
$\supp \mu$ is determined (via an explicit formula) by combinatorics (more
precisely, by the itineraries of $p$ and of the critical point).
\end{enumerate}

\end{thm}

\comm{
\begin{rem}

The statistical analysis of \cite {AM1}
and \cite {AM4} goes much beyond the Collet-Eckmann property, and
establishes a very detailed description of
maps in $\SSS \cap \R$.
In particular, we get the following additional properties of the dynamics of
$f_c$ for almost every $c \in \R$ which
is not regular or infinitely renormalizable,
\begin{enumerate}
\item The critical point is polynomially recurrent with exponent $1$:
$$
\limsup \frac {-\ln |f_c^n(0)|} {\ln n}=1,
$$
\item The critical orbit is equidistributed with respect to the absolutely
continuous invariant measure $\mu$:
$$
\lim \frac {1} {n} \sum_{i=0}^{n-1}
\phi(f_\lambda^i(0))=\int \phi d\mu
$$
for any continuous function $\phi:I \to \R$,
\item The Lyapunov exponent of the critical value
$$
\lim \frac {1} {n} \ln |Df^n(f(0))|,
$$
exists and coincides with the Lyapunov exponent of $\mu$.
\item The Lyapunov exponent of a periodic point $p$ contained in
$\supp \mu$ is determined (via an explicit formula) by combinatorics (more
precisely, by the itineraries of $p$ and of the critical point).
\end{enumerate}

\end{rem}
}

\subsection{Zero area (Theorem \ref {measure zero})}

Again, by renormalization considerations, Theorem \ref {measure zero}
reduces to the statement that (in a Multibrot-like family) almost every
non-renormalizable parameter is regular.  In view of Corollary \ref {aesss}, it
is thus enough to prove the following statement:

\begin{thm} \label {SSSzero}

The set $\SSS$  has zero area.

\end{thm}

\begin{pf}

Fix an arbitrary $c_0 \in \SSS $ and let  $\VV^{n+1}=\VV^{n+1}(c_0)$.
For $c\in \VV^{n+1}$, let $g_{n,c}$ denote
the first return map to $V^n$ under iteration by $f_c$,
let $V^n_*(c)$ be the component of $D_n(c)=\Dom (g_{n,c})$  
containing the critical value $g_{n,\tl c}(0)$.%
\footnote{In what follows we let $D_n=D_n(c_0)$, and use the similar
convention for other objects moving over $\VV^{n+1}$.}
Let 
$$
   \VV^{n+1}_* = \{c\in \VV^{n+1}:\, g_{n,c} (0)\in V^n_*(c)\}
$$
By the ``Density Points Argument'' of the previous section,
it is sufficient to show that 
\be\label{density}
   \limsup \P(\SSS |\, \VV^{n+1}_*)<1 .
\ee   

For $c\in \VV^{n+1}$, let $Z^{n+1}(c)$ be the union of the
boundaries of puzzle pieces that are subordinate to depth (and hence move
holomorphically over $\VV^{n+1}$)
over $\VV^{n+1}$.                                 
Persistent puzzle pieces include all components of $ D_n(c)$.
By Slodkovski's Theorem,
the holomorphic motion of $Z^{n+1}$ extends to a holomorphic motion $h$ 
of the whole complex plane $\C$. 

The map $\psi:c \mapsto g_{n,c} (0)$
is a diagonal to the tubing of $\partial V^n_*$ over $\VV^{n+1}_*$, so $h$
and $\psi$ give rise to a phase-parameter map
$\chi_n : V^n_* \to \VV^{n+1}_*$.  
Since 
$$
   \mod (V^n (c)\sm V^n_*(c))\geq \de n
$$ 
for all $c \in \VV^{n+1}$, the first statement of Lemma \ref {1}
implies that $\mod(\VV^{n+1} \setminus
\VV^{n+1}_*) \to \infty$,  so by the $\la$-lemma,  
$\chi_n$ is $\gamma_n$-qc, where $\gamma_n\to 1$.%
\footnote{This kind of rules relating the dynamical and parameter objects are described in \cite{AM3}
  as the {\it Phase-parameter relation}. For most purposes, one can use use these rules axiomatically.}

Given two measurable  sets $X\subset Y$ and a bi-measurable injection $\phi: Y\ra \C$,
we let 
$$
    \P_\phi(X\,|\,Y)=\P(\phi(X)\,|\,\phi(Y))
$$
be the $\phi$-pullback of the  conditional probability.  
Let $\gamma>1$.
Given a Jordan disk $V$ and a measurable set $X\subset D$,
let us define the {\it $\gamma$-capacity} $\P_\gamma(X\,|\,V)$ as follows:
$$
\P_\gamma(X\,|\,V)=\sup \P_\phi(X\,|\, V),
$$
where $\phi$ ranges over all $\gamma$-quasiconformal homeomorphisms 
$V \to \phi(V) \ssubset \C$.
Clearly, the $\gamma$-capacity is a conformal invariant. 
Let $\alpha_n = \P_\gamma(D_n\,|\, V^n)$.

For $n>0$, the set $D_n$ is ``uniformly porous'' in $V^n$  in the following sense:
There  exist $K>0$,  $\mu>0$ and $\eta > 0 $ such that
any component $W$ of $D_n$ is contained in the nest of two topological disks,   
$W\subset W' \subset W''\subset V^n$, such that:

\begin{itemize}

\item $\mod(W''\sm W')\geq \mu$;

\item $ W'$ is a $K$-quasidisk;
   
\item  $\P(D_n|\, W')\leq 1-\eta$.

\end{itemize}
To obtain such a nest,  take the return map $g_n=g_{n,c}: W\ra V^n$,
extend it to a branched covering $\hat g_n:  W''\ra V^{n-1}$ of degree $d$ or 1,
and let $W'$ be the pullback by $\hat g_n$ 
of a  big  intermediate quasidisk $U$,  $V^n\Subset U\Subset V^{n-1}$.
Since $D_{n-1}$ is not dense in $U$ (once $U$ is sufficiently big),
the Koebe Distortion Theorem implies that $D^n$ has a gap of a definite size in $W'$
(compare Lemma B.3 of \cite {AM3}). 

Uniform porosity implies that $\alpha_n<1$ for $n>0$ 
(making use of the Besikovic Covering Lemma). 

Since $g_{n,c}: V^n_* (c)\ra V^n(c)$ is a conformal map for $c\in \VV^{n+1}$,
the connected components of the set
$$
      \De_n(c)= (g_{n,c}|V^n_*(c))^{-1} (D_n(c))
$$
are puzzle pieces which are subordinate to depth $n+1$
and hence $\De_n(c)$ is respected by the holomorphic motion $h$.
Moreover, for $c \in \SSS\cap \VV^{n+1}$, $g_{n,c}(0)\in
\De_n(c)$.
It follows that $\SSS \cap \VV^{n+1}_*$ is contained in the image of
 $\De_n$ under the phase-parameter map $\chi_n :V^n_* \to \VV^{n+1}_*$.
Since this map is $\gamma$-quasiconformal for large $n$,  
$\P(\SSS|\VV^{n+1}_*) \leq \alpha_n$ by definition of the capacity.  
Thus, to prove (\ref{density}), it is enough to show that 
\be\label{alpha}
          \limsup \alpha_n<1.
\ee
We will obtain this by means of the following simple statistical argument.

 Let $\Om^{n+1}= (g_n|\, V^{n+1})^{-1}(D_n)$.  
For each connected component $W$ of $\Om^{n+1}$, 
we have: 
\be\label{delta n}
\mod(V^{n+1} \setminus W) \geq \frac{1}{d^2} \mod(V^{n-1} \setminus V^n) \geq \delta n. 
\ee

Call a component of $\Om^{n+1}$ critical if it contains $0$ and precritical
if its image under $g_n$ contains $0$.  Let $E^{n+1}$ be the union of
critical and precritical components.
If $s_n=0$ (the central return case) then $E^{n+1}=V^{n+2}$; 
otherwise $E^{n+1}$ consists of $d+1$ puzzle pieces.  
In any case, $E^{n+1}$ is the union of at most $d+1$ puzzle pieces 
$W_i \subset V^{n+1}$, 
each satisfying (\ref{delta n}).  It follows that
\be
   \P_\gamma(E^{n+1}|\, V^{n+1}) 
   \equiv  \epsilon_{n+1} \leq e^{-\delta n}.
\ee

Furthermore, 
if $W$ is a connected component of $\Om^{n+1} \setminus E^{n+1}$ 
then $g_n^2 : W \to V^n$ is a conformal map, 
and $g_n^2 (W \cap D_{n+1}) \subset  D_n$.
It follows that if $\phi:V^{n+1} \to \phi(V^{n+1}) \ssubset \C$ is a
$\gamma$-qc homeomorphism then for any such component $W$ we
have $\P_\phi(D_{n+1}|W) \leq \alpha_n$ (by the definition of capacity).
Hence
$$
\P_\phi(D^{n+1} \,|\,\Omega^{n+1} \setminus E^{n+1}) \leq \alpha_n,
$$
so that,
$$
   \P_\phi ( D^{n+1}  \setminus E^{n+1}  |\, V^{n+1}) \leq
(1-\P_\phi(E^{n+1}\,|\,V^{n+1})) \alpha_n.
$$
Thus,
\begin{align*}
 \P_\phi ( D_{n+1}\,|\, V^{n+1} )&=  \P_\phi(E^{n+1}\,|\,V^{n+1})
   +\P_\phi(D_{n+1} \setminus E^{n+1}\,|\,V^{n+1})\\
\nonumber
&\leq \P_\phi(E^{n+1}\,|\,V^{n+1})+(1-\P_\phi(E^{n+1}\,|\,V^{n+1})
\alpha_n\\
\nonumber
&=\alpha_n+(1-\alpha_n) \P_\phi(E^{n+1}\,|\,V^{n+1}).
\end{align*}
 Taking the supremum over all $\phi$ under consideration, 
we obtain:
$$
  \alpha_{n+1} \leq \alpha_n+(1-\alpha_n) \epsilon_{n+1},
$$
 so 
$$
      \frac {1-\alpha_{n+1}} {1-\alpha_n} \geq 1-\epsilon_{n+1} \geq 1-e^{-
\de 
n},
$$
which yields (\ref{alpha}).
\end{pf}

\subsection{Porousity of $\MM$ (Theorem \ref {refined})}

If $c$ is not combinatorially recurrent, 
then by Theorem \ref{Misiur}  $c$ is a Lebesgue density point of the complement of $\MM$.

So, assume that $c$ is combinatorially recurrent.  By Theorem \ref{bounds},
$\cap \YY^n(c)=\{c\}$, and by Corollary \ref {SSS} and Theorem \ref {SSSzero},
$\liminf \P(\MM|\YY^n(c))<1$.  This is not enough, though,
to conclude that $c$ is not a Lebesgue density point of $\MM$, since the
$\YY^n(c)$ do not in general have a bounded shape (where a
set $K \subset \C$ is said to have $C$-bounded shape if it contains a round
disk of radius $\frac {1} {C} \diam(K)$).\footnote
{One can show that the parapuzzle pieces corresponding to the
puzzle pieces in the {\it enhanced nest} constructed in \cite {KSS} have a bounded shape,
but this nest is less convenient for the statistical arguments.}  
However,  the following lemma will allow us to replace them with
shrinking  domains of bounded shape.

\begin{lemma}

For every $\delta>0$, there exists $\kappa>1$ with the
following property.  Let $D$ be a Jordan disk and let $x \in D$.  Then there
exists $r>0$ such that $\D_r(x) \subset D$ and for every Jordan disk $\De \subset D$ 
that intersects both $\D_r(x)$ and $\di \D_{\kappa r}(x)$, 
we have: $\mod(D \setminus \De) <\de$.

\end{lemma}

\begin{pf}

Let $\psi:(D,x) \to (\D,0)$ be the Riemann map, 
and let $r$ be maximal radius such that $\psi(\D_r(x)) \subset \D_{1/2}$.

Notice that if $B \subset \D$ is a Jordan disk with 
$\mod(\D \setminus B) \geq \delta$ and $B \cap \D_{1/2} \neq \emptyset$,
then $B \subset \D_R$ where $R=R(\delta)<1$.%
\footnote{The optimal choice is to take $R\in (1/2, 1) $ such
that $\mod(\D \setminus [1/2,R])=\delta$.}
By the Koebe Distortion Theorem, $\psi^{-1}(\D_R)
\subset \D_{\kappa r}(x)$.  If $\De $ intersects $ \D_r(x)$ and
$\partial \D_{\kappa r}(x)$, then $\psi(\De)$ intersects $ \D_{1/2}$
and $\partial \D_R$, so that,  
$\mod(D \setminus \De)=\mod(\D \setminus \psi(\De))<\delta$.
\end{pf}

Given $\rho>0$, there
exists $n$ such that $\VV^n \subset \D_\rho(c)$,
$\mod(\VV^n \setminus \VV^{n+1})>\delta$ and $\mod(V^n
\setminus V^{n+1})>\delta$.  Then for some $\eta=\eta(\delta)$, we
have $\mod(\VV^{n+1} \setminus \VV^{n+2}(\tilde c))>\eta$ 
and $\mod(V^{n+1}(\tilde c) \setminus V^{n+2}(\tilde c))>\eta$ 
for every combinatorially recurrent parameter $\tilde c \in \VV^{n+1}$. 
 Almost every $\tilde c \in \VV^{n+1} \cap \MM$ is combinatorially recurrent, 
and by Corollary \ref {SSS} and Theorem \ref {SSSzero},  
$\P(\MM|\VV^{n+2}(\tilde c))<1-\epsilon$.  
By the previous lemma, there exists $r>0$ such that $\D_r(c) \subset \VV^{n+1}$ and
any $\VV^{n+2}(\tilde c)$ intersecting $\D_r(c)$ is contained $\D_{\kappa r}(c)$
for some $\kappa>1$.
Let $X \subset \D_{\kappa r}(c)$ be the union of all the $\VV^{n+2}(\tilde c)$
intersecting $\D_r(c)$.  Then 
$$1-\P(\MM\, |\, \D_{\kappa r}(c)) \geq \kappa^{-2} (1-\P(\MM\, |\, X)) \geq \epsilon \kappa^{-2}.$$

\begin{rem} \label {gen}

Let us indicate how to generalize Theorem \ref {refined} to the finitely
renormalizable case.  We can not just argue via renormalization since it
would only prove that parameters are not density points of a copy
of the Multibrot set containing it, and indeed
a neighborhood of a satellite renormalizable parameter (with repelling
periodic orbits) contains
non-renormalizable parameters belonging to infinitely many Misiurewicz
limbs.

This can be solved by constructing a different version of the puzzle and
parapuzzle, which is designed to be compatible with a fixed
renormalization.  Namely, one constructs ``adapted Yoccoz puzzle pieces
of depth $0$'', where instead of using the external rays landing at the
$\alpha$-fixed point of $f_c$,
one uses the external rays landing at the orbit of the $\alpha$-fixed point
of the renormalization of $f_c$.  Though the combinatorial description is
different (see \cite {Sch2}, \cite {M2}),
the whole geometric and statistical
analysis can be carried out to obtain Theorem \ref {refined} in the more
general setting.

\end{rem}

%% file: main.bbl
\begin{thebibliography}{BKNS}

\bibitem[ALM]{ALM} A. Avila, M. Lyubich, and W. de Melo,
Regular or stochastic dynamics in real analytic families of unimodal
maps.
Invent. Math. 154 (2003), 451--550.

\bibitem[AM1]{AM1} A. Avila, C. G. Moreira.
Statistical properties of unimodal maps: the quadratic family.
Ann. Math. 161 (2005), 831-881.


\bibitem[AM3]{AM3} A. Avila, C. G. Moreira.
Phase-parameter relation and sharp statistical properties for general
families of unimodal maps.
Contemp. Math. 389 (2005), 1-42.

\bibitem[AM4]{AM4} A. Avila, C. G. Moreira.
Statistical properties of unimodal maps: physical measures, periodic orbits
and pathological laminations.
Pub. Math. IH\'ES 101 (2005), 1-67.

\bibitem[AKLS]{AKLS} A. Avila, J. Kahn, M. Lyubich, and W. Shen,
Combinatorial rigidity for unicritical polynomials.  Preprint IMS at Stony
Brook, \# 5 (2005).  To appear in Annals of Math.


\bibitem[BC]{BC} M. Benedicks, L. Carleson. On iterations of $1-ax^2$
 on (-1,1). Ann. Math. 122 (1985), 1-25.

\bibitem[BR]{BR} L. Bers and H. L. Royden,  Holomorphic families of
injections.  Acta Math.  157  (1986),  no. 3-4, 259--286.

\bibitem[BKNS]{BKNS}
H. Bruin, G. Keller, T. Nowicki, S. van Strien.
Wild Cantor attractors exist. Ann. of Math. (2) 143 (1996), 97--130.

\bibitem[BSS1]{BSS1}
H. Bruin, W. Shen, S. van Strien.
Invariant measures exist without a growth condition.  Comm. Math. Phys.  241
(2003),  no. 2-3, 287--306.

\bibitem[BSS2]{BSS}
H. Bruin, W. Shen, S. van Strien.
Existence of unique SRB-measures is typical for real unicritical polynomial
families.  Ann. Sci. \'Ecole Norm. Sup. (4)  39  (2006),  no. 3, 381--414.

\bibitem[DH1]{DH} Douady and Hubbard.
Etude dynamiques des polyn\^omes complexes. 
Publications Math. d'Orsay, 84-02 (1984) and 85-04 (1985). 

\bibitem[DH2]{DH-pl} Douady and Hubbard.
  On the dynamics of polynomial-like maps. 
  Ann. Sci. \'Ec. Norm.  Sup., v. 18 (1985), 287--344. 

\bibitem[GM]{GM} Goldberg, Lisa R.; Milnor, John
Fixed points of polynomial maps. II. Fixed point portraits.
Ann. Sci. \'Ecole Norm. Sup. (4)  26  (1993),  no. 1, 51--98.

\bibitem[H]{H} J. H. Hubbard, Local connectivity of Julia sets and
bifurcation loci: three theorems of J.-C. Yoccoz.  Topological methods in
modern mathematics (Stony Brook, NY, 1991),  467--511, Publish or Perish,
Houston, TX, 1993.

\bibitem[J]{J} M. Jakobson. Absolutely continuous invariant measures
for one-parameter families of one-dimensional maps. Comm. Math. Phys.,
 v. 81 (1981), 39-88.

\bibitem[K]{K} J. Kahn. Holomorphic removability of Julia sets.
 Preprint IMS at Stony Brook, \# 11 (1998). 

\bibitem[KL1]{KL1} J. Kahn and M. Lyubich,
The quasi-additivity law in conformal geometry.  Preprint IMS Stony Brook
\#2 (2005).  To appear in Annals of Math.

\bibitem[KL2]{KL2} J. Kahn and M. Lyubich,
Local connectivity of Julia sets for unicritical polynomials.
Preprint IMS Stony Brook \#3 (2005).  To appear in Annals of Math.

\bibitem[KSS]{KSS} O. Kozlovski, W. Shen, and S.~van~Strien. 
  Rigidity for real polynomials.   Ann. of Math. (2)  165  (2007),  no. 3,
749--841.

\bibitem[L1]{L7} M. Lyubich.  On the Lebesgue measure of the Julia set
of a quadratic polynomial.  Preprint IMS at Stony Brook, \# 10 (1991).

\bibitem[L2]{attractors} M. Lyubich. Combinatorics, geometry and attractors
of quasi-quadratic maps. Ann. Math, {\bf 140} (1994), 347-404.

\bibitem[L3]{puzzle} M. Lyubich,
Dynamics of quadratic polynomials. I, II.
Acta Math.  178  (1997),  no. 2, 185--247, 247--297.

\bibitem[L4]{parapuzzle} M. Lyubich. Dynamics of quadratic polynomials, III.
Parapuzzle and SBR measure.
Ast\'erisque,
v. 261 (2000),  173 - 200.

\bibitem[L5]{horseshoe} M. Lyubich. Almost every real quadratic maps is either regular 
 or stochastic. Ann. Math., v. 156 (2002), 1--78. 

\bibitem[MSS]{MSS} R. Ma\~n\'e, P. Sad, and D. Sullivan.
 On the dynamics of rational maps, Ann. scient. Ec. Norm. Sup.,
  {\bf 16}  (1983), 193-217.

\bibitem[M1]{M-attractors} J. Milnor. On the concept of attractor. 
  Comm. Math. Physics., v. 99 (1985), 177--195. 

\bibitem[M2]{M1} J. Milnor. Local connectivity of Julia sets: Expository lectures.
In: ``The Mandelbrot set: themes and variations. 
London Math. Soc. Lecture Notes, v. 274 (ed. Tan Lei), p. 67-116. 

\bibitem[M3]{M2} J. Milnor. Periodic orbits,  external rays, and the Mandelbrot set:
   Expository account. In: G\'eometrie Complexe et Systems Dynamiques.
   Ast\'erisque, v. 261 (2000), 277--333 . 

\bibitem[MN]{MN} M. Martens and T. Nowicki. 
  Invariant measures for Lebesgue typical quadratic maps.
Ast\'erisque, v. 261 (2000),  239--252. 

\bibitem[RL]{RL} J. Rivera-Letelier, On the
   continuity of Hausdorff dimension of Julia sets and similarity between
   the Mandelbrot set and Julia sets. Fund. Math. 170 (2001), no. 3,
   287--317.


\bibitem[Sa]{Sa} D. Sands.  Misiurewicz maps are rare.
   Comm. Math. Phys. 197 (1998), no. 1, 109--129.

\bibitem[Sc1]{Sch1} D. Schleicher. Rational parameter rays of the Mandelbrot set.
Ast\'erisque,
v. 261 (2000), 405-443.

\bibitem[Sc2]{Sch2} D. Schleicher. On fibers and renormalization of Julia sets 
  and Multibrot sets.  Preprint IMS at Stony Brook, \# 13 (1998).

\bibitem[S]{S} M. Shishikura.  Topological, geometric and
complex analytic
properties of Julia sets.  In Proceedings of the International Congress
of Mathematicians (Z\"urich, 1994), pp. 886-895.  Birkh\"auser, Basel
(1995).

\bibitem[Sl]{Sl} Z. Slodkowsky. Holomorphic motions and polynomial hulls.
Proc. Amer. Math. Soc.,  111 (1991), 347--355.

\bibitem[Sm]{Sm} D. Smania.
Puzzle geometry and rigidity: The Fibonacci cycle is hyperbolic.
J. Amer. Math. Soc.  20  (2007),  no. 3, 629--673

\bibitem[vS]{vS} S. van Strien. Misiurewicz points unfold generically
  (even if they are critically non-finite). Fundamenta Math. 163 (2000),
   39--54. 


\bibitem[T]{T} M. Tsujii.  Positive Lyapunov exponents in families of one
dimensional dynamical systems.  Invent. Math. 111 (1993), 113-137.

\end{thebibliography}
